# Circle Packing Problem Using Nature-Inspired Optimization Techniques

Veni Goyal[a*], Pulkit Mundra[a*], Kusum Deep[b]

[a] Mechanical and Industrial Engineering Department, Indian Institute of Technology Roorkee, Roorkee – 247667, Uttarakhand, India

[b] Department of Mathematics, Indian Institute of Technology Roorkee, Roorkee – 247667, Uttarakhand, India

*Dual First Authorship

**Abstract -** This paper deals with the problem of circle packing, in which the largest radii circle is to be fit in a confined space filled with arbitrary circles of different radii and centers. A circle packing problem is one of a variety of cutting and packing problems. We suggest four different nature-inspired Meta-heuristic algorithms to solve this problem. Algorithms are based on the social behavior of other biology species such as birds, wolves, fireflies, and bats. Moreover, recent advancements in these algorithms are also considered for problem-solving. The circle packing problem is one of the NP-hard problems. It is challenging to solve NP-hard problems exactly, so the proposed algorithms provide an approximate solution within the allotted time. Standard statistical parameters are used for comparison, and simulation and results indicate that the problem is highly non-linear and sensitive.

## 1. Introduction

### 1.1. Nature-Inspired Optimization Techniques

These algorithms are highly efficient in finding optimized solutions to multi-dimensional and multi-modal problems [1]. The conventional optimization approach in calculus is finding the first-order derivative of the objective function and equating it to zero to get the critical points. These critical points then give the maximum or minimum value as per the objective function. Calculating gradients or even higher-order derivatives need more computing resources and is more error-prone than other methods.

Further, you can imagine how complex it is to find the solution to a minimization/ maximization problem with 20 or even more variables. However, using these nature-inspired algorithms can solve the problem with less computational effort and time complexity. These algorithms use a stochastic approach to find the best solution in the large search space of the problem.

### 1.2. Circle Packing problem

The problem is described as searching for the largest radii circle to be fit in a confined space filled with arbitrary circles of different radii and centers. The solution circle should entirely lie within the bounds of the confined space without overlapping already existing circles. In geometry, **circle packing** is the study of the arrangement of circles (of equal or varying sizes) on a given surface such that no overlapping occurs and so that no circle can be enlarged without creating an overlap [2].

Application – Radiation Treatment Planning: Expose the affected area so that the surrounding organs at risk are least affected.

Example: Suppose the kidney bean shape is the affected area, and we have to expose it with radiation and one stroke is in the form of a circle. Here our objective is to minimize the number of strokes, cover the affected area as much as possible and prevent radiation exposure to its surrounding as it can be fatal to health.

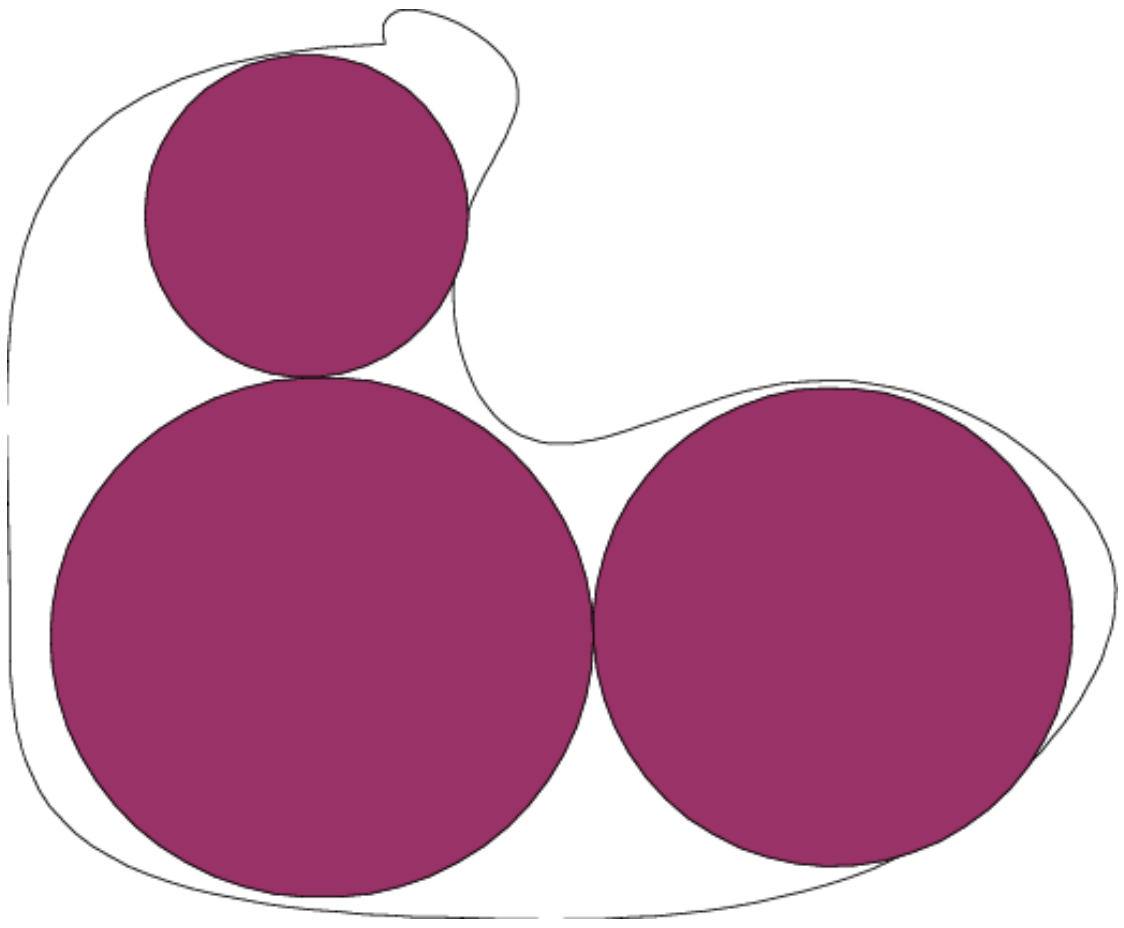

**Figure 1** Circle Packing in Kidney Bean

## 2. Mathematical Modelling of Circle Packing Problem

The mathematical model for the circle packing problem involves finding a maximum radius circle at a point A that fits in the confined space satisfying the necessary condition of no overlapping. For this, we have calculated the minimum distances between the point and all the existing circles and also the distance from the space bounds and finally taken the least value among these distances as the largest radius for the circle to be packed in the space with center at the point X.

Mathematically, let coordinates of point X as ($X_1$, $X_2$). We know the centers ($C_i$) and radii ($R_i$) of all the pre-existing circles. $D_i$ can be defined as the minimum distance from the point to the circle, i.e., the normal distance. In vector form,

$D_i$ = sqrt (sum ((X – $C_i$) ^2)) – $R_i$

As the space is confined and for simplicity, we have chosen rectangular bounds. Let $B_i$ be the distances of the point X from the upper and lower bounds.

$B_1 = X_1 - lb$
$B_2 = ub - X_1$
$B_3 = X_2 - lb$
$B_4 = ub - X_2$

Finally, we have taken the minimum distance value as the largest radius for the circle that can be packed in space with its center at point X.

$R_o$ = Minimum {$D_i$, $B_i$}

**Table 1** Abbreviations

| | |
|---|---|
| PSO | Particle Swarm Optimization |
| PSOd | Particle Swarm Optimization with the new update mechanism |
| CPSO | Constricted Particle Swarm Optimization |
| GWO | Grey Wolf Optimization |
| RWGWO | Random Walk Grey Wolf Optimization |
| ApFA | Firefly Algorithm with adaptive control parameters |
| FA | Firefly Algorithm |
| BA | Bat Algorithm |
| rand | Random Number in [0,1] |
| sqrt | Square Root |
| sum | Summation |
| ub | Upper bound |
| lb | Lower bound |
| wMax | Inertia maximum |
| wMin | Inertia minimum |
| vMax | Velocity maximum |
| vMin | Velocity minimum |

## 3. Optimization Algorithms

### 3.1. Particle Swarm Optimization

In 1995, Kennedy and Eberhart [3] presented Particle Swarm Optimization (PSO) inspired by the social behaviour of fish or birds to solve non-linear global optimization problems. PSO uses individual agents or particles to search for the solution in objective function space. Each particle has its personal best location, and all particles have their global best location. Every particle gets attracted to the individual best and global best position while maintaining the tendency to move randomly. Any particle finds a better place; then it will update its position to a better one.

New velocity vector and position factor are given by following formula:

$$v_i^{t+1} = (v_i^t + \phi_1 U_1^t (pb_i^t - x_i^t) + \phi_2 U_2^t (gb_i^t - x_i^t))$$

$$x_i^{t+1} = x_i^t + v_i^{t+1}$$

where $gb_i^t$ and $pb_i^t$ are global best and personal best locations and $\phi_1$ and $\phi_2$ are two parameters called acceleration coefficients $U_1^t$ and $U_2^t$ are two d × d diagonal matrices with diagonal elements distributed in the interval [0, 1) uniformly at random.

**Table 2** PSO parameters and values

| Parameters | Values |
|---|---|
| wMax | 0.9 |
| wMin | 0.2 |
| $c_1$ | 2.0 |
| $c_2$ | 2.0 |
| vMax | (ub – lb) *(0.2) |
| vMin | -vMax |

### 3.2. Particle Swarm Optimization with a New Update Mechanism

Since its invention, many PSO variants have been proposed by modifying its solution update rule to improve its performance. In 2017, Kiran [4] proposed a new update in PSO to prevent stagnation of the particle population and remove parameters such as $c_1$ and $c_2$ which, are social and cognitive components of the velocity update rule. PSOd proposed a new position update rule based on the normal distribution. The position of the particle is obtained by the normal distribution.

$X_{i,j}(t+1) = \mu + \sigma * Z$

where, $\mu$ (mean), $\sigma$ (standard deviation) and Z is given by:

$\mu = (X_{i,j}(t) + pbest_{i,j}(t) + gbest(t))/3$

$\sigma = sqrt\ [\ ((X_{i,j}(t) - \mu)^2 + (pbest_{i,j}(t) - \mu)^2 + (gbest_i(t) - \mu)^2)/3]$

$Z = (-2ln(k_1))^{1/2} x\ cos(2\pi k_2)$

where, $k_1$ and $k_2$ are uniform random numbers produced in range of [0,1].

**Table 3** PSOd parameters and values

| Parameters | Values |
|---|---|
| $c_1$ | random |
| $c_2$ | random |
| z | Sqrt ((-1) * log($c_1$)) * cos (2* (3.14) *$c_2$) |

**3.3 Constricted Particle Swarm Optimization**

Eberhart [5] used a constriction factor to the particle velocity to prevent unlimited growth of particle velocity and ensure the convergence of particle swarm optimizer.

Modified velocity equation:

$v_i^{t+1} = \chi\ (v_i^t + \phi_1 U_1^t\ (pb_i^t - x_i^t) + \phi_2 U_2^t\ (lb_i^t - x_i^t))$

with $\chi = 2/|2 - \varphi - sqrt\ (\varphi^2 - 4\varphi)|$ where $\chi$ is the constriction factor and $\varphi$ = = $c_1 + c_2$ and $\varphi > 4$, usually $c_1$ and $c_2$ are set to 2.05 [6]. It can be concluded that the best approach to use with particle swarm optimization as a "rule of thumb" is to utilize the constriction factor approach while limiting $V_{max}$ to $X_{max}$.

**Table 4** CPSO parameters and values

| Parameters | Values |
|---|---|
| wMax | 0.9 |
| wMin | 0.2 |
| $c_1$ | 2.05 |
| $c_2$ | 2.05 |
| $c_3$ | $c_1 + c_2$ |
| A | 2/ \| (2 – $c_3$ – sqrt($c_3$^2 – 4*$c_3$) \| |
| vMax | (ub – lb) * 0.2 |
| vMin | -vMax |

**3.4 Grey Wolf Optimization**

In 2014, Mirjalili [7] presented Grey Wolf Optimizer (GWO) based on the social hierarchy and hunting behavior of wolves. The GWO algorithm is a new contribution to the family of swarm intelligence-based metaheuristics. In the family of Swarm intelligence algorithms, GWO is the only algorithm based on position managers. In addition to wolves' social hierarchy, group hunting is another interesting social behavior of grey wolves. For designing the mathematical model, we have to consider social hierarchy. The Fittest solution will be given by alpha, the second fittest by beta, the third best by delta, and the remaining are given by the omega. Alpha *(α),* Beta (*β)*, and Delta *(δ)* will lead the hunting, and omega *(ω)* will follow them. It is assumed that alpha, beta, delta have a better knowledge of the solution, which is the prey's location, and omega will update their position according to the three best solutions.

The encircling strategy by the wolves around the prey is mathematically modelled by proposing the following equations as:

$X_{t+1} = X_{p,t} - \mu * d$

$d = |c * X_{p,t} - X_t|$

$\mu = 2 * b * r_1 - b$

$c = 2 * r_2$

where $X_{t+1}$ is the position of the wolf at $(t + 1)^{th}$ iteration $X_t$ is the position of the wolf at $t^{th}$ iteration, $X_{p,t}$ is the position of the prey at $t^{th}$ iteration, d is difference vector, μ and c are coefficient vector and b is the linearly decreasing vector from 2 to 0 over iterations, expressed as:

$b = 2 - 2.\ (t/\text{max number of iterations})$

The hunting strategy of the grey wolves can be mathematically modeled by approximating the prey position with the help of α, β, and δ and solutions (wolves). Therefore, by following this approximation, each wolf can update their positions by:

$X_1' = X_\alpha - \mu_\alpha * d_\alpha$

$X_2' = X_\beta - \mu_\beta * d_\beta$

$X_3' = X_\delta - \mu_\delta * d_\delta$

$X_{t+1} = (X_1' + X_2' + X_3')/3$

where $X_\alpha$, $X_\beta$, $X_\delta$ are the positions approximated by α, β and δ solutions.

**3.5 Random Walk Grey Wolf Optimization**

In GWO proposed by mirjalili, the wolves' positions are updated by the leading wolves, but how will alpha update its position as it is the leader of the pack. Due to this, there is premature convergence at local optima, and GWO cannot converge to global optima. In 2019, Gupta and Deep [8] presented a novel Random Walk Grey Wolf Optimizer to solve this problem to update wolves' positions differently. It updates the position of alpha, beta, and delta by random walk, and for omega, it updates its position as the given method in GWO.

Random walk is a random process that consists of consecutive random steps. The relationship between any two consecutive random walks can be defined as:

$W_N = \sum^{N}_{i=1} (s_i) = W_{N-1} + S_N$

where $s_i$ is a random step taken from any random distribution, this relationship shows that the next state $W_N$ is only dependent on the current state $W_{N-1}$ and the step taken from the current state to the next state. The step size can be fixed or can vary. So, for a wolf starting with a point $x_0$ and suppose its final location is $x_N$, then a random walk can also be defined as:

$x_N = x_0 + \alpha_1 s_1 + \alpha_2 s_2 + \ldots\ldots\ldots + \alpha_N s_N = x_0 + \sum^{N}_{i=1}(\alpha_i s_i)$

where $\alpha_i>0$ is a parameter that controls the step size $s_i$ in each iteration.

### 3.6 Firefly Algorithm

The Firefly Algorithm (FA) is a swarm intelligence-based optimization technique developed by Yang [9]. FA work on three principles: 1. Fireflies are unisex, which means any firefly could be attracted to another firefly regardless of their sex; 2. Attractiveness is proportional to their brightness; 3. The prospect of the objective function determines the light intensity of a firefly.

The flashing light can be formulated so that it is associated with the objective function to be optimized. In the maximization type of problem, the brightness can be related to the value of the objective function.

The distance between any two fireflies i and j at $x_i$ and $x_j$, respectively, is the Cartesian distance:

$R_{ij} = ||x_i - x_j|| = sqrt\ (\sum^{d}_{k=1}(x_{ik} - x_{jk})\ \string^2)$

The movement of a firefly i is attracted to another more attractive (brighter) firefly j is determined by:

$x_i = x_i + \beta_0 e^{-\gamma r_{ij}^2} (x_j - x_i) + \alpha(\varepsilon_i)$

where the second term is due to the attraction while the third term is randomization with α being the randomization parameter α ∈ [0, 1] and $\varepsilon_i$ is a random value uniformly distributed in the range [−0.5, 0.5].

**Table 5** FA parameters and values

| Parameters | Values |
|---|---|
| alpha | 1.0 |
| beta | 1.0 |
| gamma | 0.1 |
| betamin | 0.1 |
| theta | 0.97 |
| d | 2.0 |

### 3.7 Firefly Algorithm with adaptive control parameters

Firefly Algorithm with adaptive control parameters (ApFA) presented by Wang [10] includes some changes in the control parameters of FA. It dynamically adjusts $\alpha$ with a new adaptive parameter strategy, wherein the original FA is fixed to 1. It also alters the attractiveness coefficient $\beta_0$ with a simple dynamic approach.

Updated $\alpha$ and $\beta_0$ are as follows:

*$\alpha (t + 1) = (1- (1/G_{max}))\, \alpha(t)$*

where $\alpha (0) = 0.5$ is used.

$\beta_0(t+1) = rand1\ (rand2<0.5)$

$\beta_0(t+1) = \beta_0(t)$ *otherwise*

where rand1 and rand2 are two random numbers generated by the uniform distribution, and the initial $\beta_0(0) = 1.0.$

**Table 6** ApFA parameters and values

| Parameters | Values |
|---|---|
| alpha | 0.5 |
| beta | 1.0 |
| gamma | 0.1 |
| betamin | 0.1 |
| theta | 0.97 |
| d | 2.0 |

### 3.8 Bat Algorithm

Bat Algorithm (BA), based on bats' echolocation behavior, is proposed by Yang [11]. BA work on three principles: 1. As mentioned, Bat uses echolocation to sense distance. They can differentiate their prey and analyze background barriers somehow; 2. Bats fly randomly with velocity $v_i$ at position $x_i$ with a fixed frequency $f_{min}$, varying wavelength $\lambda$, and loudness $A_0$ to search for prey. They can automatically adjust the wavelength (or frequency) of their emitted pulses and adjust the rate of pulse emission $r \in [0,1]$, depending on the proximity of their target; 3. Although the loudness can vary in many ways, we assume that the loudness varies from a significant (positive) $A_0$ to a minimum constant value $A_{min}$.

The new solutions $x_i^t$ and velocities $v_i^t$ at time step t is given by:

$f_i = f_{min} + (f_{max} - f_{min}) * \beta,$

$v_i^t = v_i^{t-1} + (x_i^t - x_*) * f_i$

$x_i^t = x_i^{t-1} + v_i^t$

where $\beta \in [0,1]$ is a random vector drawn from a uniform distribution. Here $x_*$ is the current global best location (solution) after comparing all the solutions among all the n bats.

**Table 7** BA parameters and values

| Parameters | Values |
|---|---|
| A | 1.0 |
| $r_0$ | 1.0 |
| alpha | 0.97 |
| gamma | 0.1 |
| $f_{max}$ | 2.0 |
| $f_{min}$ | 0 |
| t | 0 |
| d | 2.0 |

## 4. Experiment

Nature-inspired algorithms are among the most powerful algorithms for optimization. We have solved this problem using different Nature-Inspired Optimization Techniques (NIOT) such as PSO, CPSO, PSOd, GWO, RWGWO, FA, ApFA, and BA. For experimentation purpose we have used the rectangular test data set as shown in Table 8 and Table 9. Upper bound and lower bound specifies the search space for our problem while the existing circles and randomly distributed over the search space.

**Table 8** Common Parameters

| Parameters | Value |
|---|---|
| Upper bound | 100 |
| Lower bound | -100 |
| No of variable | 2 |

**Table 9** Test Data Set

| S. No. | Centre | Radius |
|---|---|---|
| 1 | [0,0] | 15 |
| 2 | [-50 0] | 12 |
| 3 | [-70 30] | 15 |
| 4 | [40 -70] | 10 |
| 5 | [20 30] | 20 |
| 6 | [60 60] | 20 |
| 7 | [50 0] | 15 |
| 8 | [-70 -30] | 15 |
| 9 | [-40 70] | 20 |
| 10 | [-20 -30] | 5 |

We experimented with six combinations of Maximum Iterations (100, 500, 1000) and Number of Particles (50, 100), and correspondingly, for each combination, we ran each algorithm for 100 different initial seeds that are randomly generated and reported Best, Worst, Mean, Median and Standard Deviation values.

**Figure 2** Test data plotted on graph

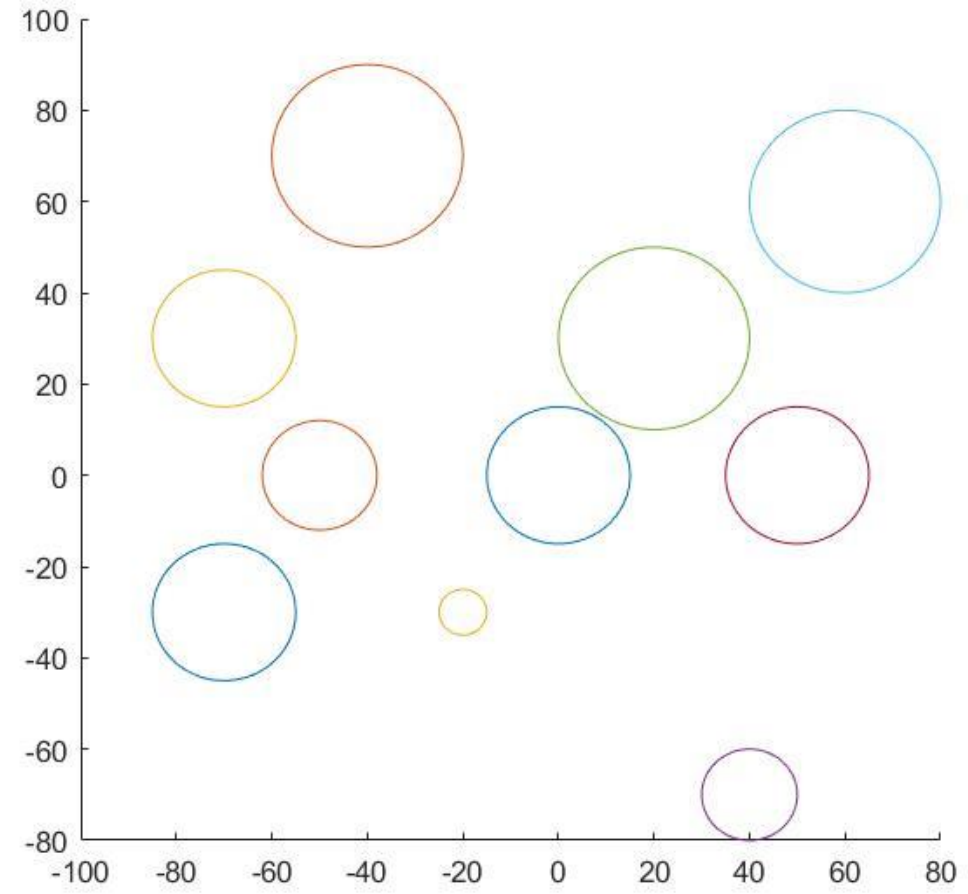

For a better comparison, we pre-generated the 100 seeds for each combination. We then imported them into all the algorithms so that the initial seed for each algorithm remains the same. The algorithms proposed in this paper are implemented in MATLAB language. These algorithms have been run on the MATLAB R2021a. Specifications of the server is Operating System: Windows 10 Home Edition, System type:64-bit Operating system, Processor: Intel® Core™ i5-7300HQ, Clock Speed :2.5Ghz, RAM:8GB.

## 5. Results and Discussions

The results of all combinations are summarised in the tables, followed by scatter plots and performance comparison graphs.

We have further discussed the behavior and performance of algorithms.

**Table10** Max Iterations = 100 & Number of Particles = 50

| | PSO | PSOd | CPSO | GWO | RWGWO | FA | ApFA | BA |
|---|---|---|---|---|---|---|---|---|
| Max Iterations | 100 | 100 | 100 | 100 | 100 | 100 | 100 | 100 |
| Number of Particles | 50 | 50 | 50 | 50 | 50 | 50 | 50 | 50 |
| Best Value | **34.2393** | **34.2393** | **34.2393** | 34.2348 | 34.235 | 25.3691 | **34.2393** | **34.2393** |
| Worst Value | 34.2056 | 30.776 | 34.2056 | 33.9978 | 30.745 | **19.9309** | 34.2056 | 22.932 |
| Mean | **34.2277** | 34.1302 | 34.2218 | 34.2004 | 34.15884 | 22.2332 | 34.22312 | 33.08165 |
| Median | **34.2393** | 34.2106 | 34.2056 | 34.20065 | 34.2004 | 22.2306 | **34.2393** | 34.2054 |
| Standard Deviation | **0.01599** | 0.37187 | 0.01692 | 0.029697 | 0.349892 | 1.246299 | 0.016921 | 2.097384 |

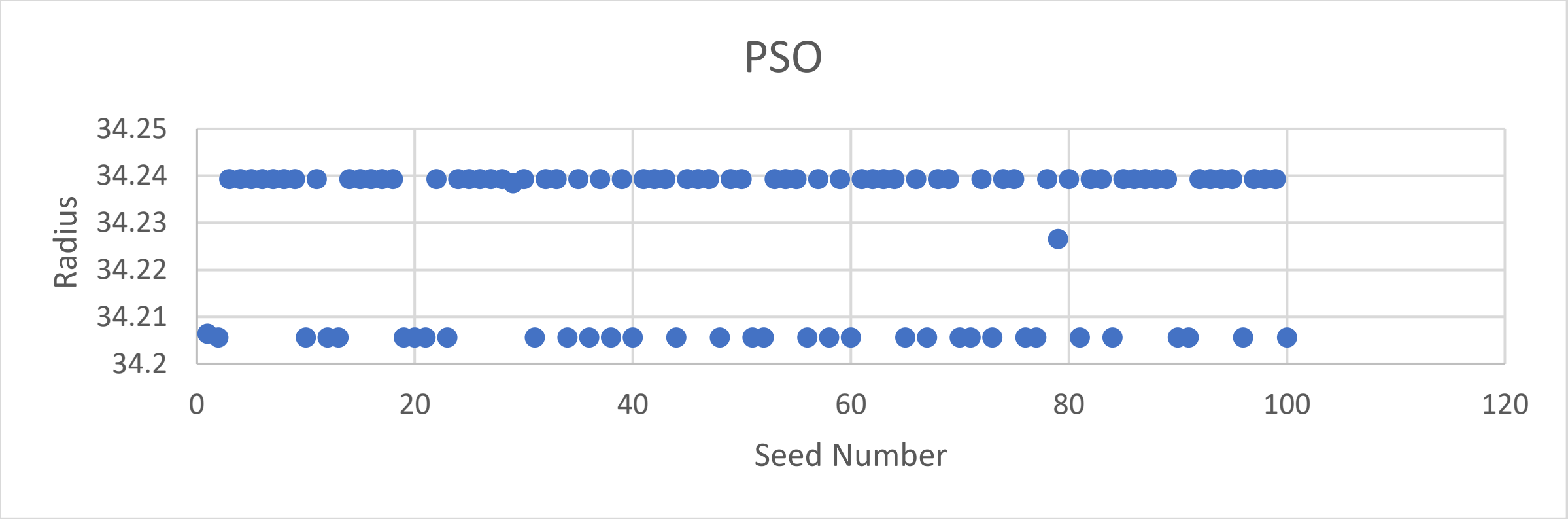

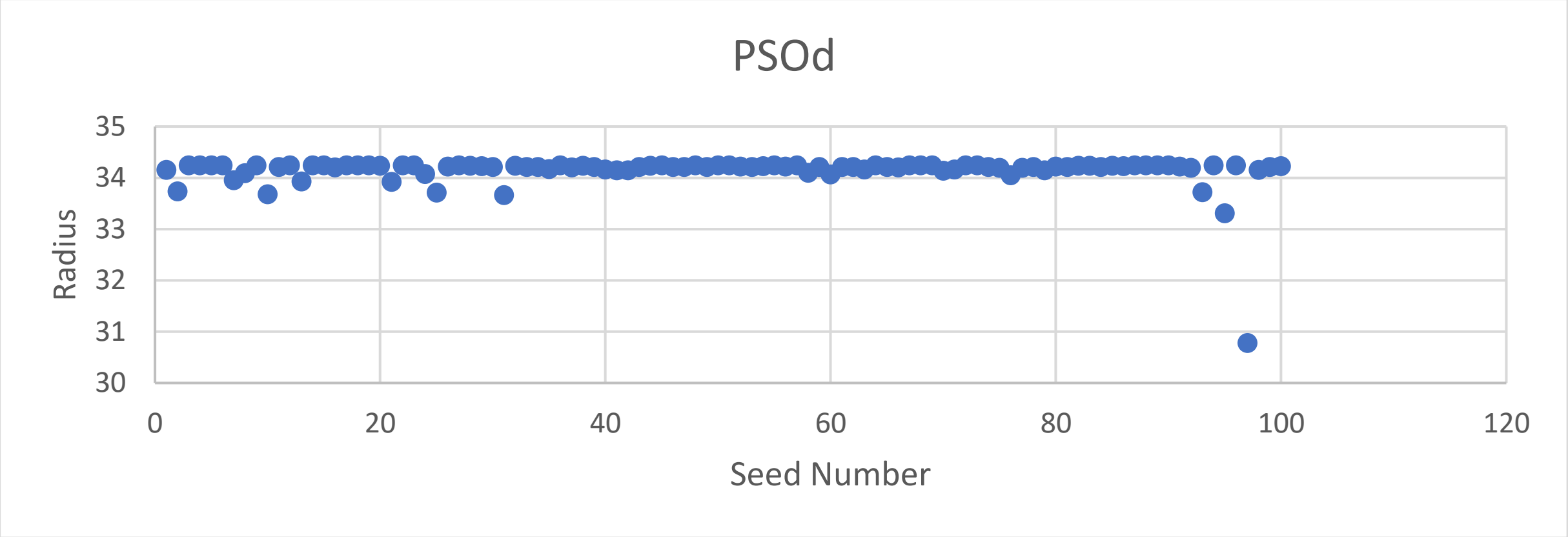

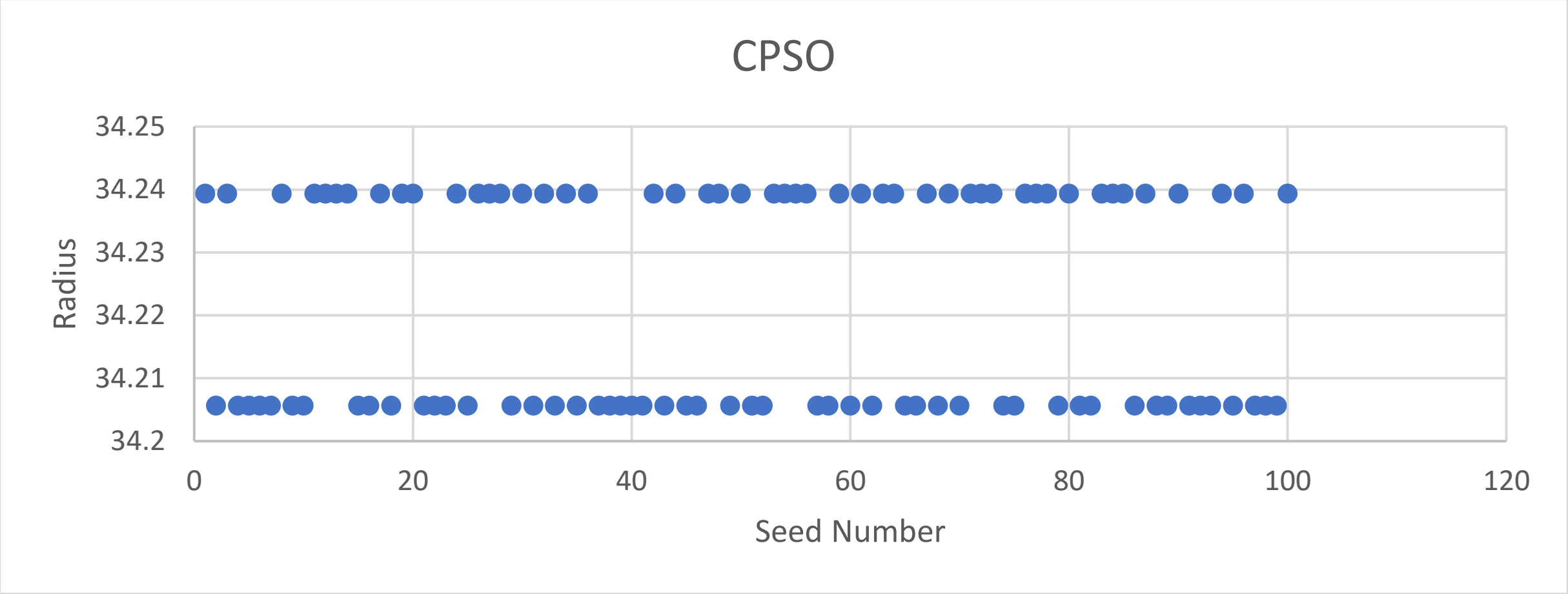

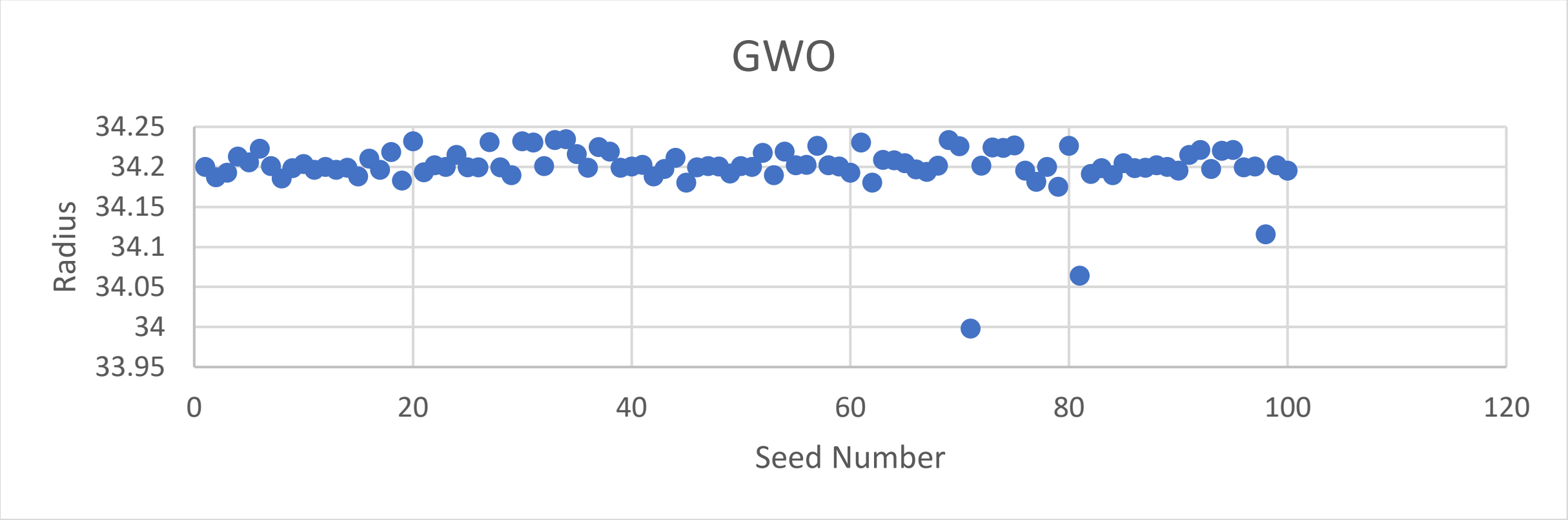

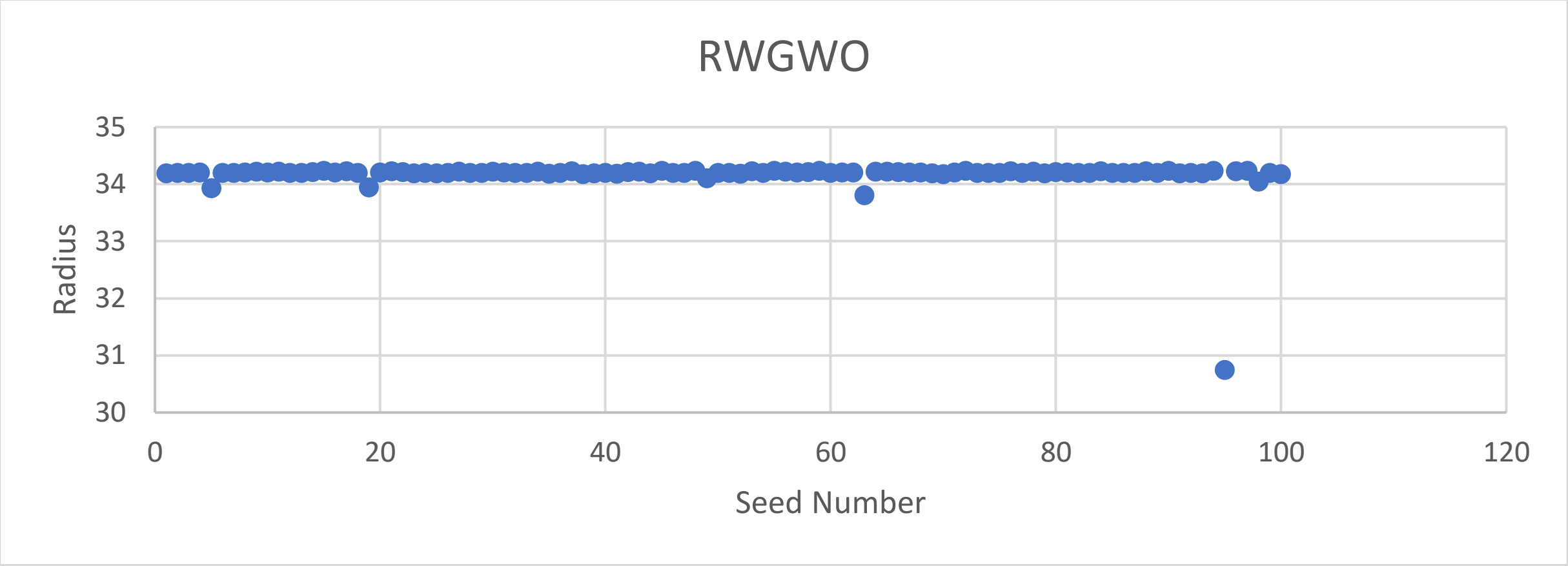

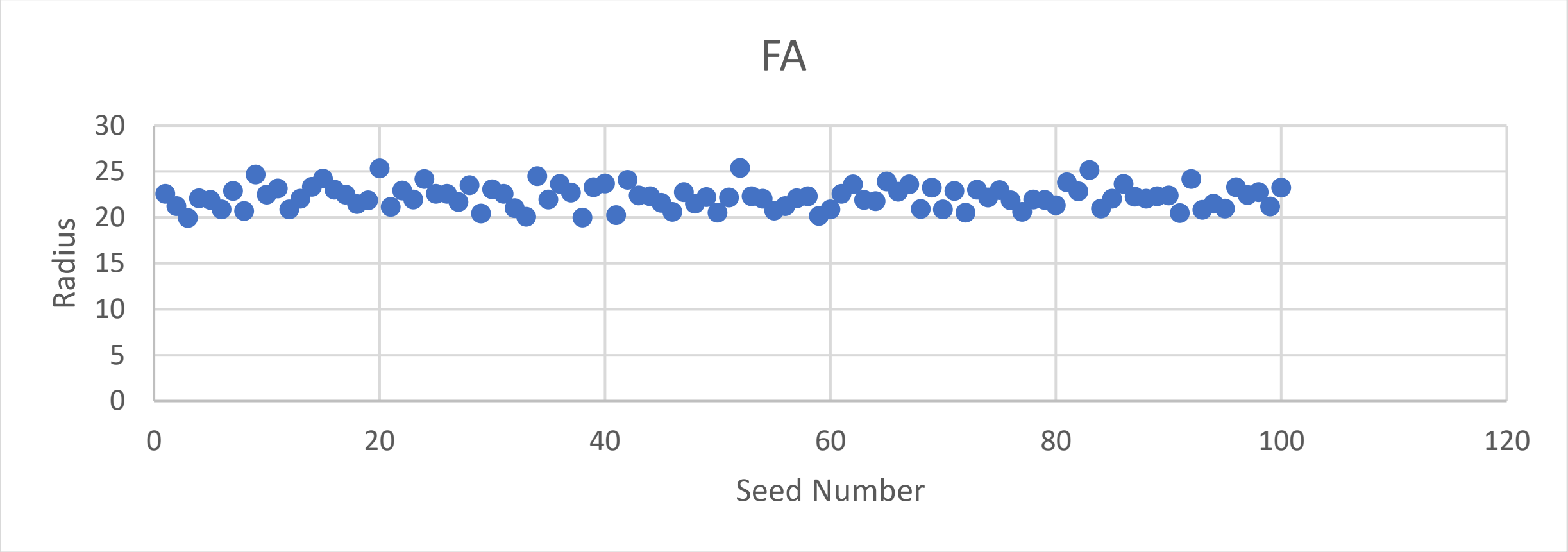

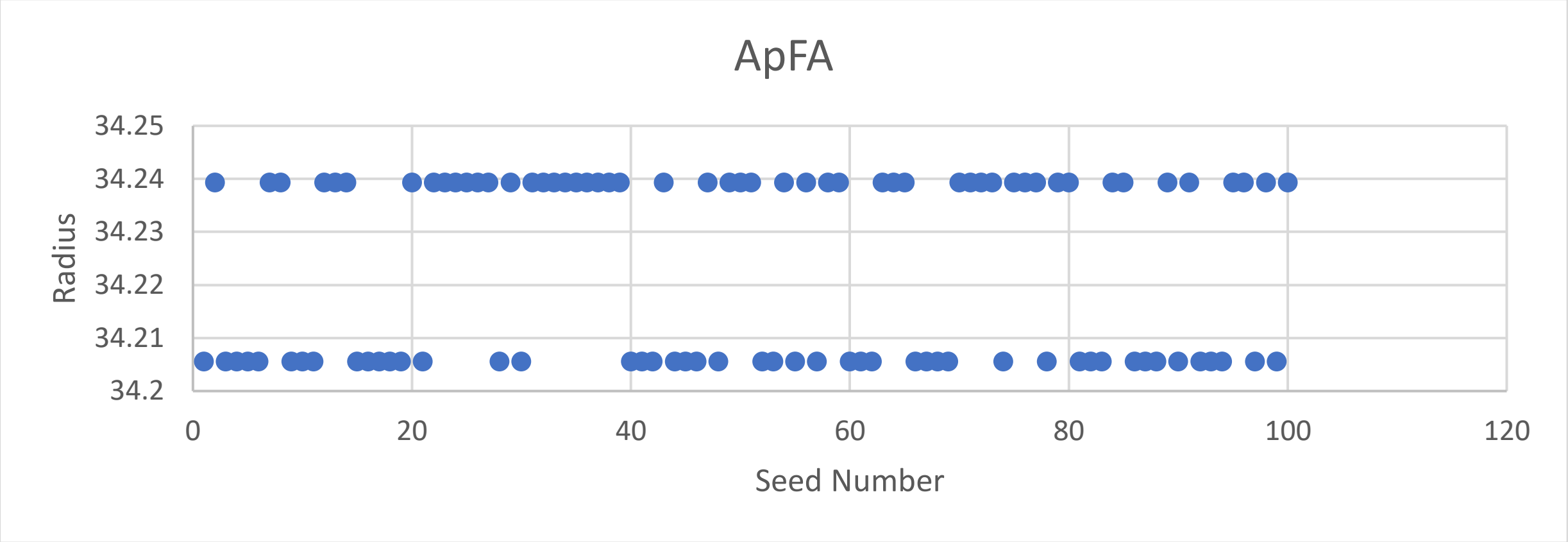

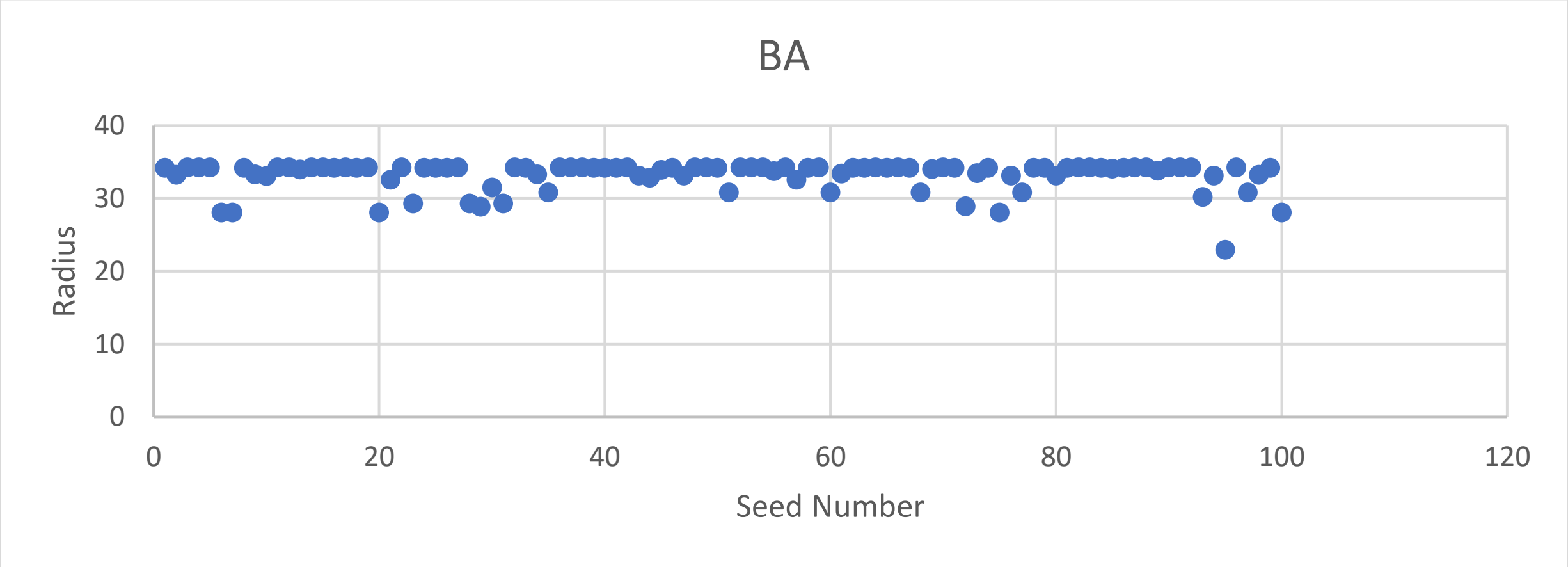

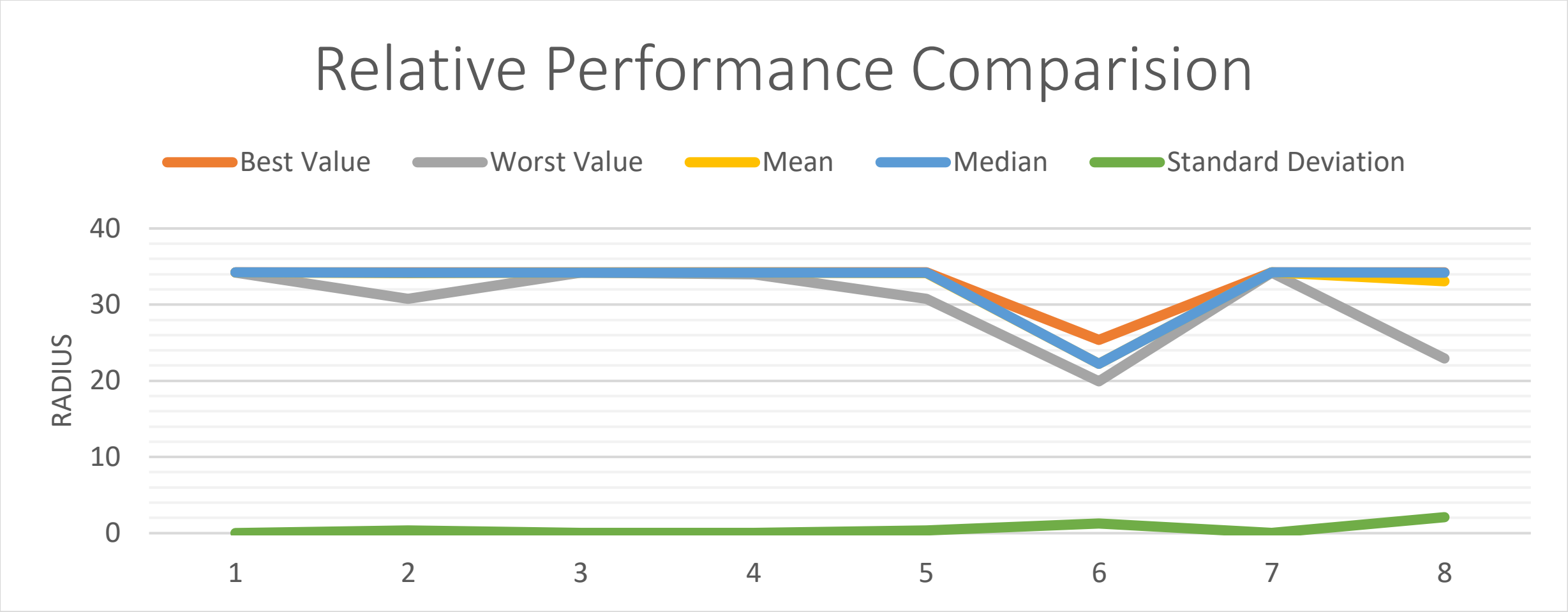

The table shows the best value for the radius is achieved by PSO, PSOd, CPSO, ApFA, and BAT algorithms. FA algorithms achieve the worst value. In terms of average and standard deviation, although there is a very minimal difference between PSO, CPSO, and ApFA, PSO has been the best-fit algorithm for solving the circle packing problem with these parameters. Also, poor results of the Firefly algorithm can be attributed to the fact that it converges slowly. Lastly, with the most significant value for standard deviation BAT algorithm has quiet instability.

Table11: Max Iterations = 500 & Number of Particles = 50

| | PSO | PSOd | CPSO | GWO | RWGWO | FA | ApFA | BA |
|---|---|---|---|---|---|---|---|---|
| Max Iterations | 500 | 500 | 500 | 500 | 500 | 500 | 500 | 500 |
| Number of Particles | 50 | 50 | 50 | 50 | 50 | 50 | 50 | 50 |
| Best Value | **34.2393** | **34.2393** | **34.2393** | 34.2384 | 34.2387 | **34.2393** | **34.2393** | **34.2393** |
| Worst Value | 30.7809 | 30.7805 | 34.2056 | 34.0772 | 34.0625 | 34.2055 | 34.2056 | **23.3138** |
| Mean | 34.1987 | 34.0833 | 34.2275 | 34.21081 | 34.20986 | **34.23045** | 34.22953 | 33.12929 |
| Median | **34.2393** | 34.2056 | **34.2393** | 34.2045 | 34.2045 | 34.2392 | **34.2393** | 34.2056 |
| Standard Deviation | 0.34547 | 0.52066 | 0.01615 | 0.021941 | 0.028269 | **0.014834** | 0.015369 | 2.121396 |

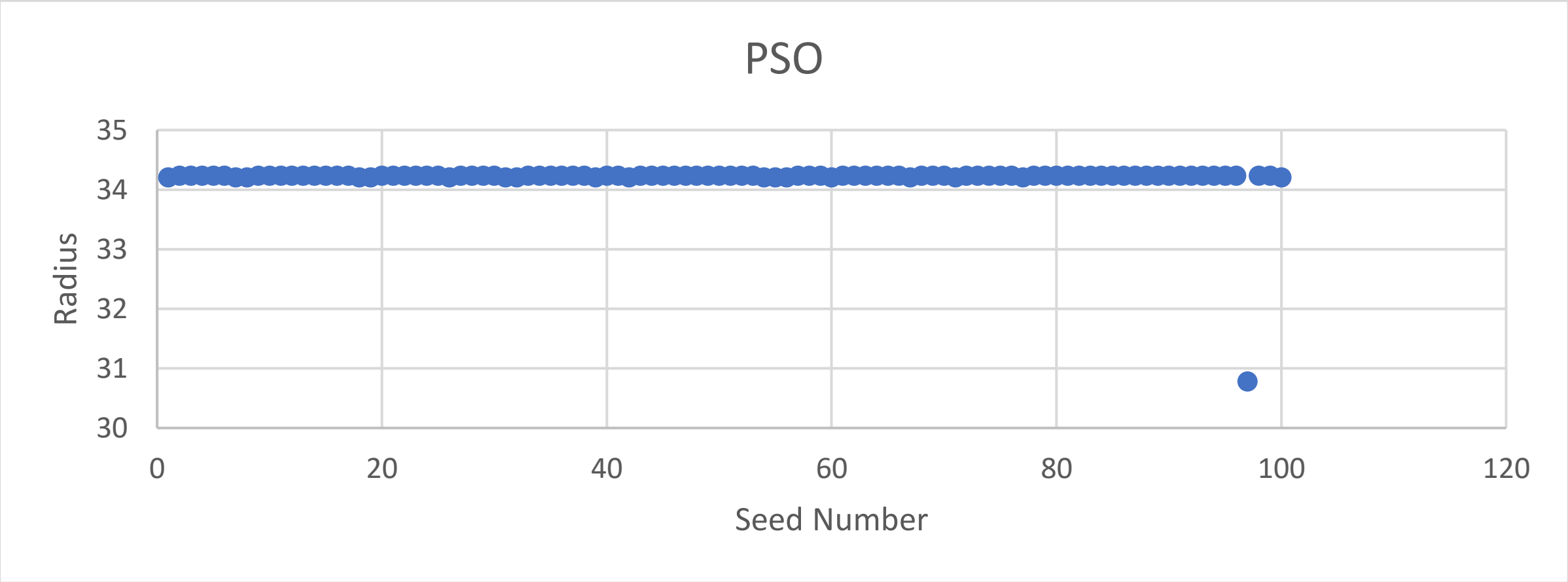

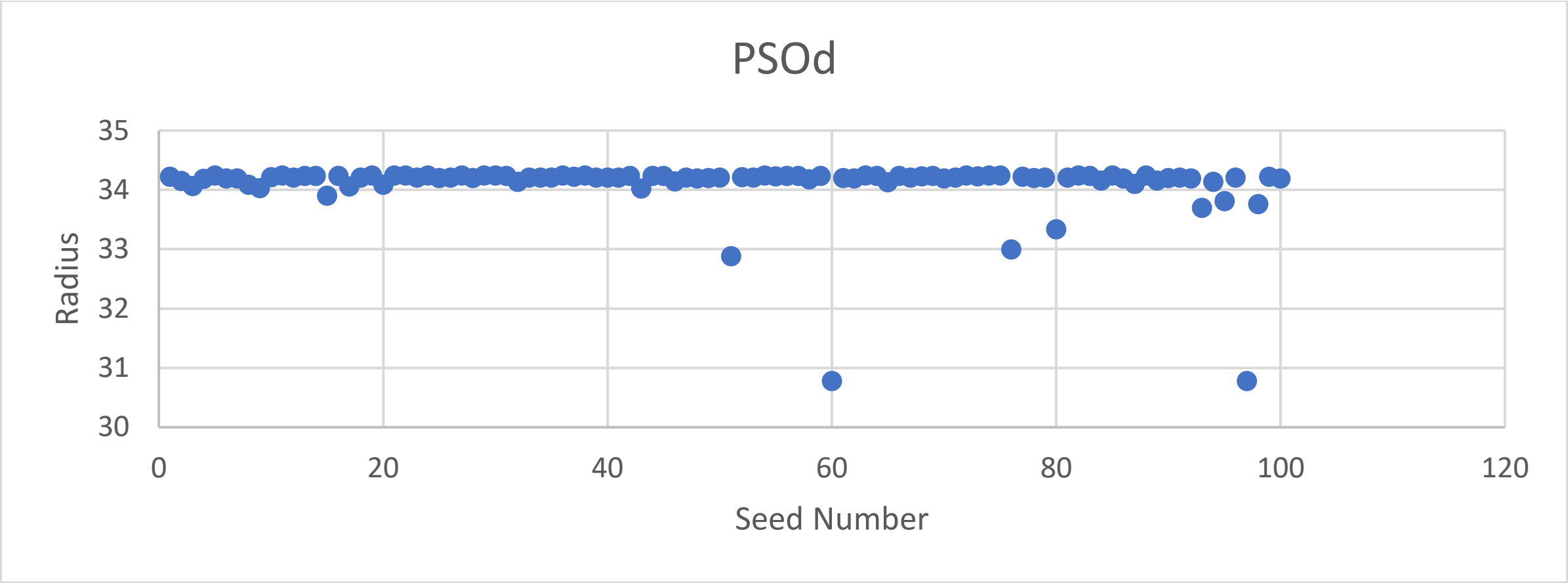
PSOd
Radius
35
34
33
32
31
30
0
20
40
60
80
100
120
Seed Number

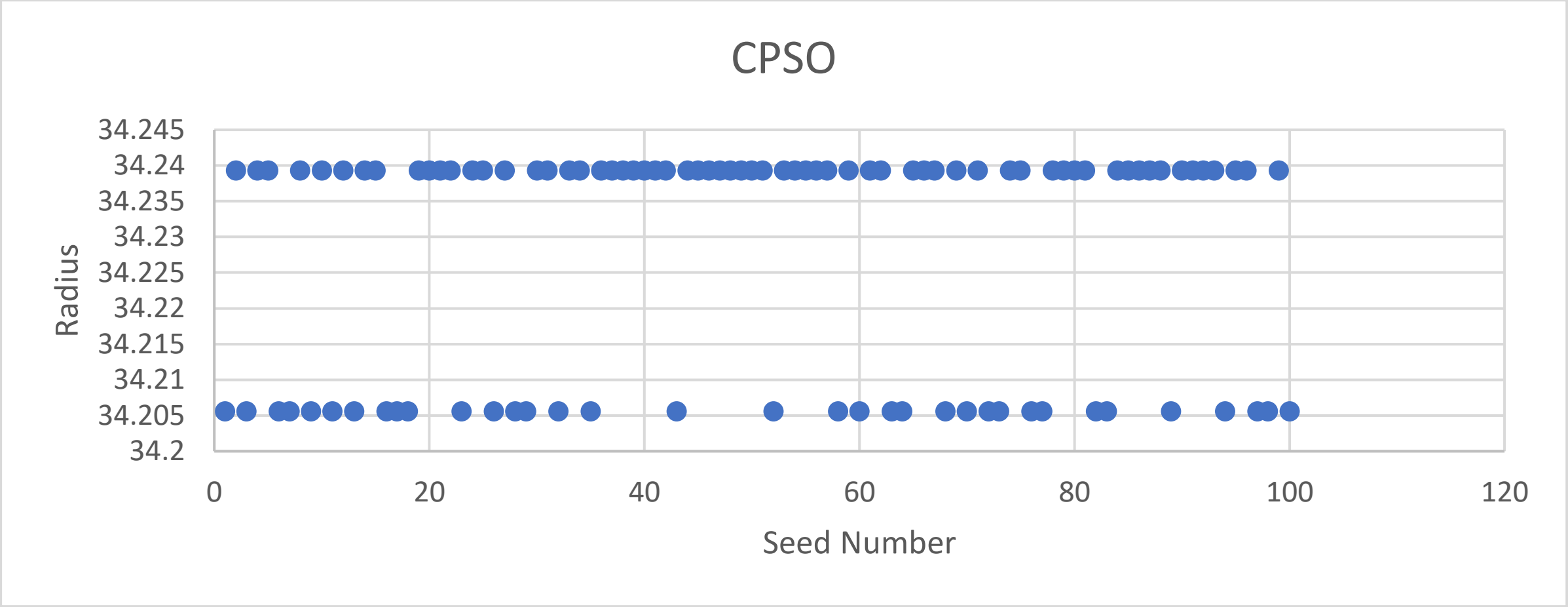
CPSO
Radius
34.245
34.24
34.235
34.23
34.225
34.22
34.215
34.21
34.205
34.2
0
20
40
60
80
100
120
Seed Number

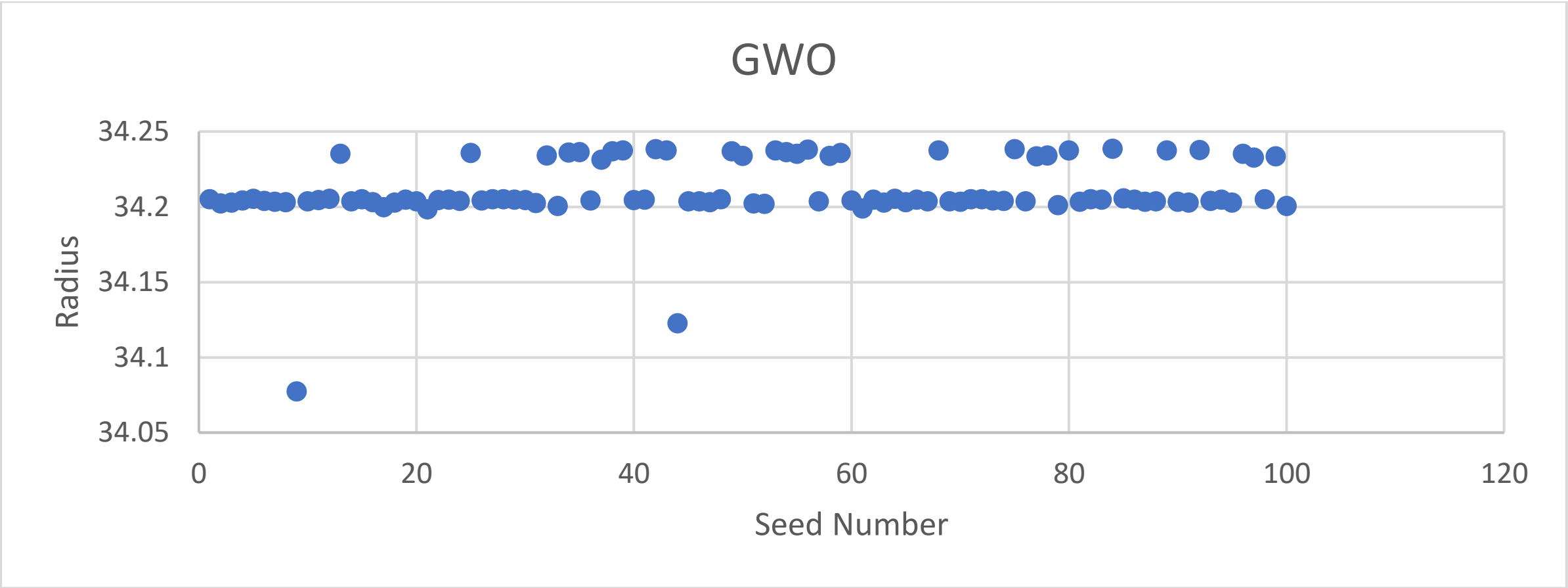
GWO
Radius
34.25
34.2
34.15
34.1
34.05
0
20
40
60
80
100
120
Seed Number

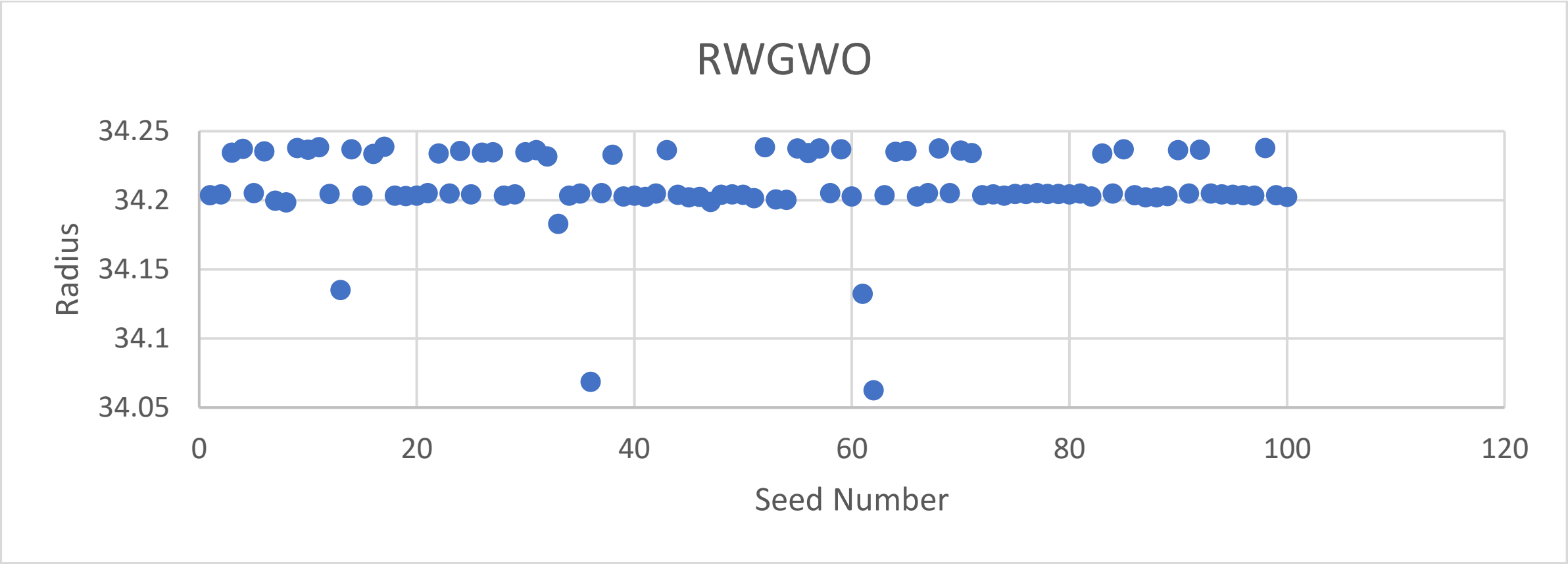
RWGWO
Radius
34.25
34.2
34.15
34.1
34.05
0
20
40
60
80
100
120
Seed Number

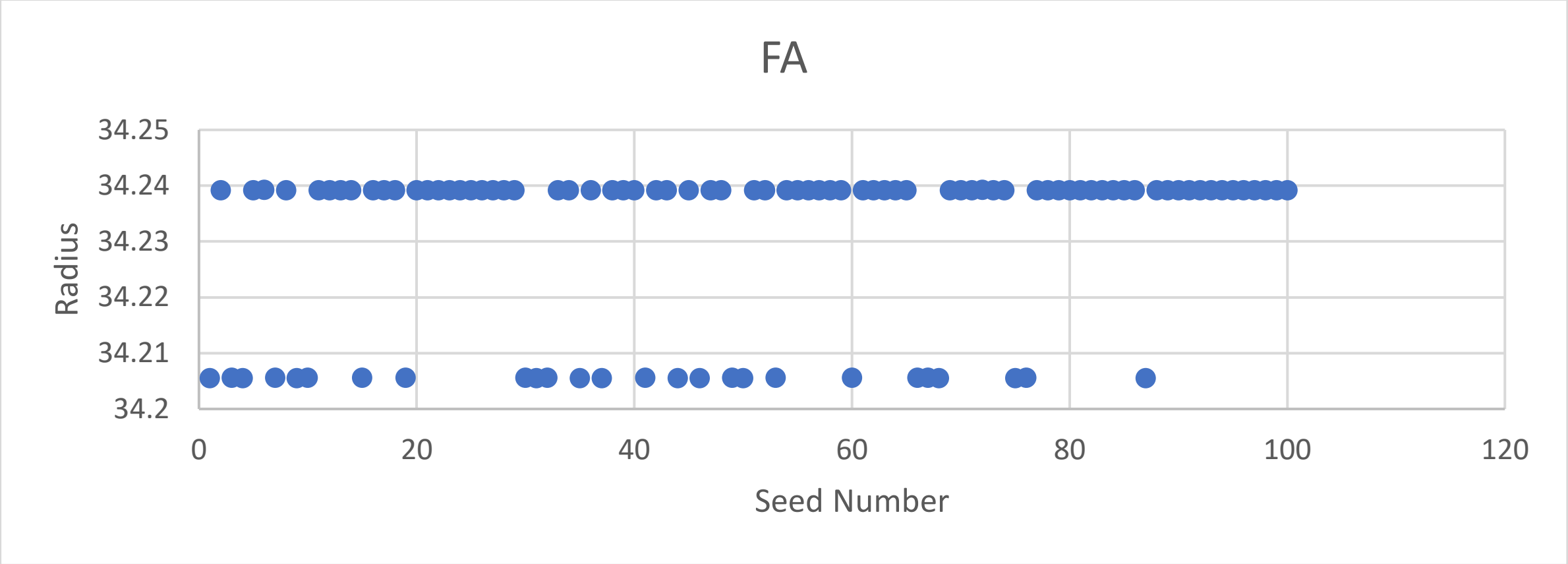
FA
Radius
34.25
34.24
34.23
34.22
34.21
34.2
0
20
40
60
80
100
120
Seed Number

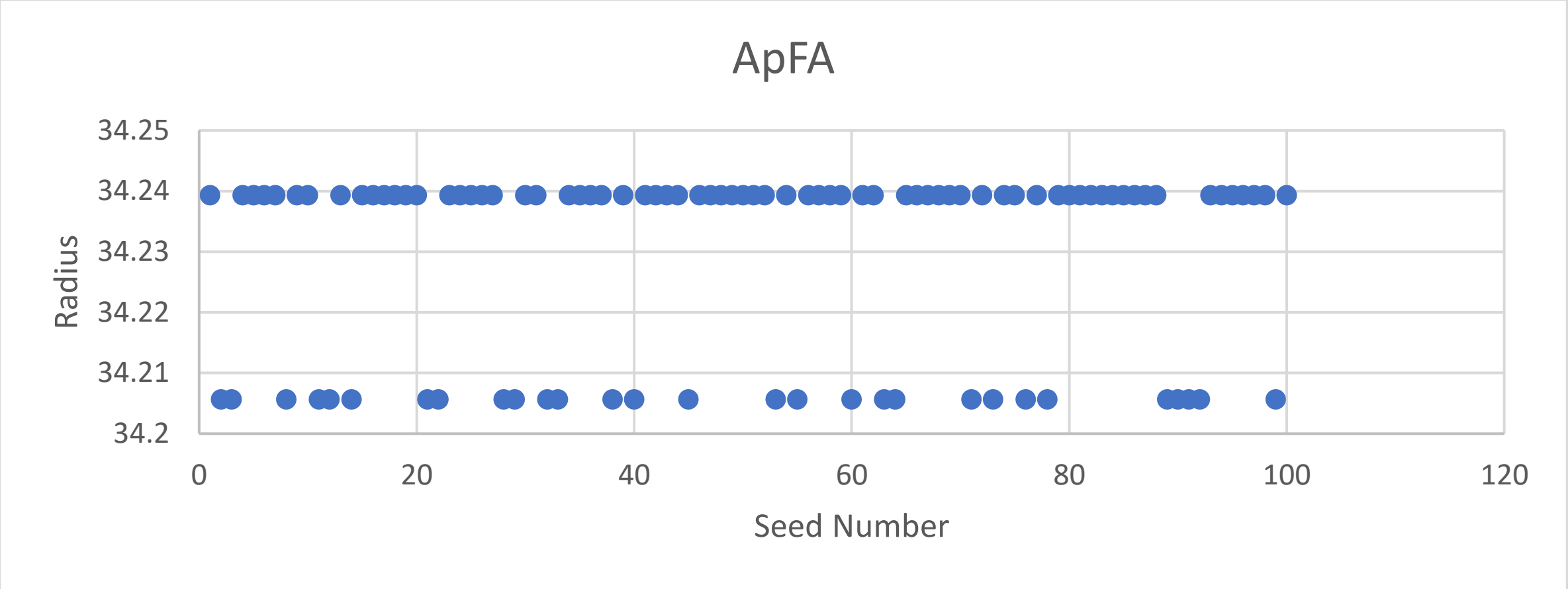
ApFA
Radius
34.25
34.24
34.23
34.22
34.21
34.2
0
20
40
60
80
100
120
Seed Number

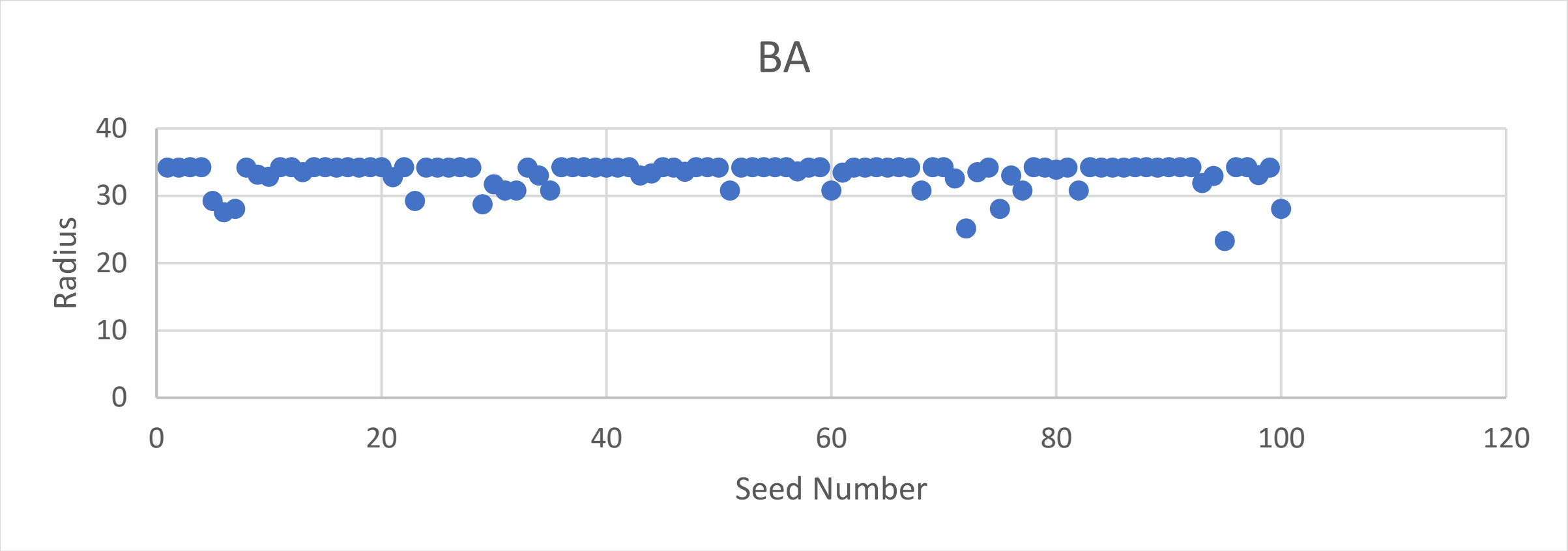
BA
Radius
40
30
20
10
0
0
20
40
60
80
100
120
Seed Number

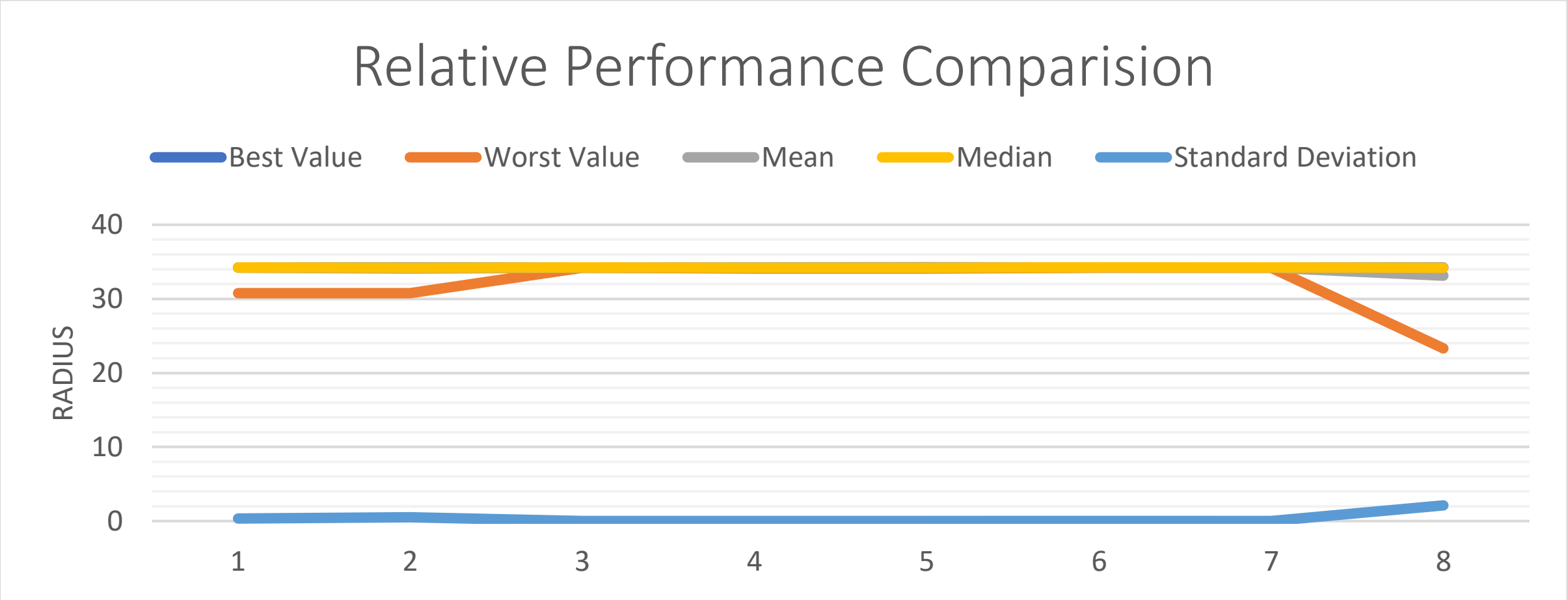

Here, in addition to PSO, PSOd, CPSO, ApFA, and BA, the best value is achieved by FA. BA achieves the worst value. In terms of average and standard deviation, there is a very minimal difference between ApFA and FA. Still, in contrast to the previous one, FA has been the best-fit algorithm for solving circle packing problems with these parameters. Notably, the performance of GWO and RWGWO has been improved. Lastly, here also with the most significant value for standard deviation BAT algorithm has depicted quiet instability.

Table12: Max Iterations = 1000 & Number of Particles = 50

| | PSO | PSOd | CPSO | GWO | RWGWO | FA | ApFA | BA |
|---|---|---|---|---|---|---|---|---|
| Max Iterations | 1000 | 1000 | 1000 | 1000 | 1000 | 1000 | 1000 | 1000 |
| Number of Particles | 50 | 50 | 50 | 50 | 50 | 50 | 50 | 50 |
| Best Value | **34.2393** | **34.2393** | **34.2393** | 34.2391 | 34.239 | **34.2393** | **34.2393** | **34.2393** |
| Worst Value | 34.2056 | 30.6411 | 30.7809 | 34.1049 | 34.0925 | 34.2056 | 34.2056 | **22.9812** |
| Mean | **34.2349** | 34.0134 | 34.1906 | 34.21381 | 34.21166 | 34.22852 | 34.22751 | 33.07364 |
| Median | **34.2393** | 34.2068 | **34.2393** | 34.2051 | 34.205 | **34.2393** | **34.2393** | 34.2056 |
| Standard Deviation | **0.01134** | 0.76059 | 0.34481 | 0.018805 | 0.024541 | 0.015799 | 0.016155 | 2.137182 |

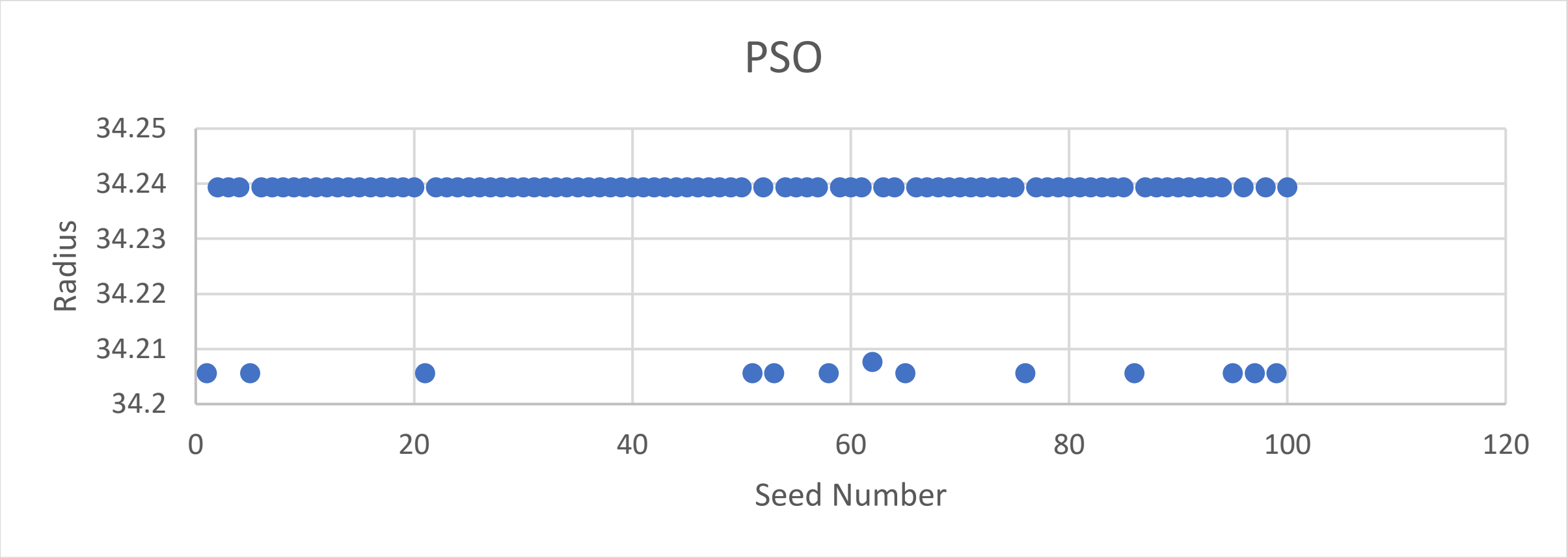

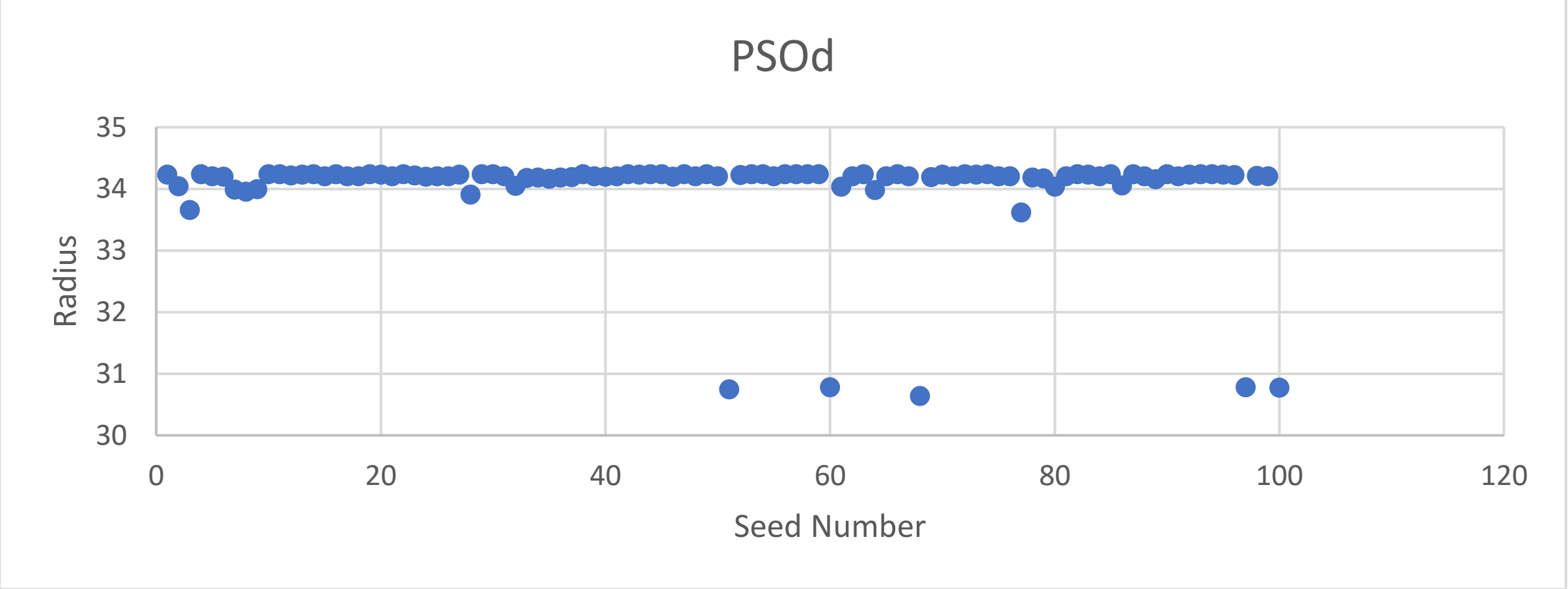

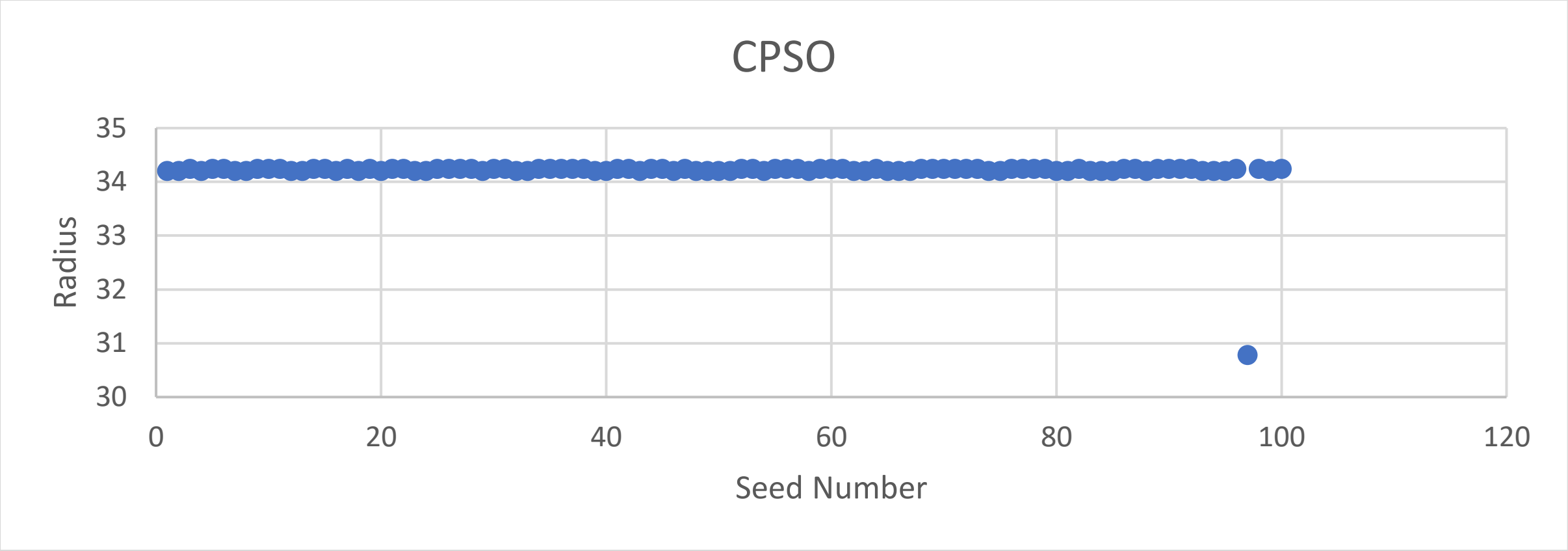
CPSO
Radius
35
34
33
32
31
30
0
20
40
60
80
100
120
Seed Number

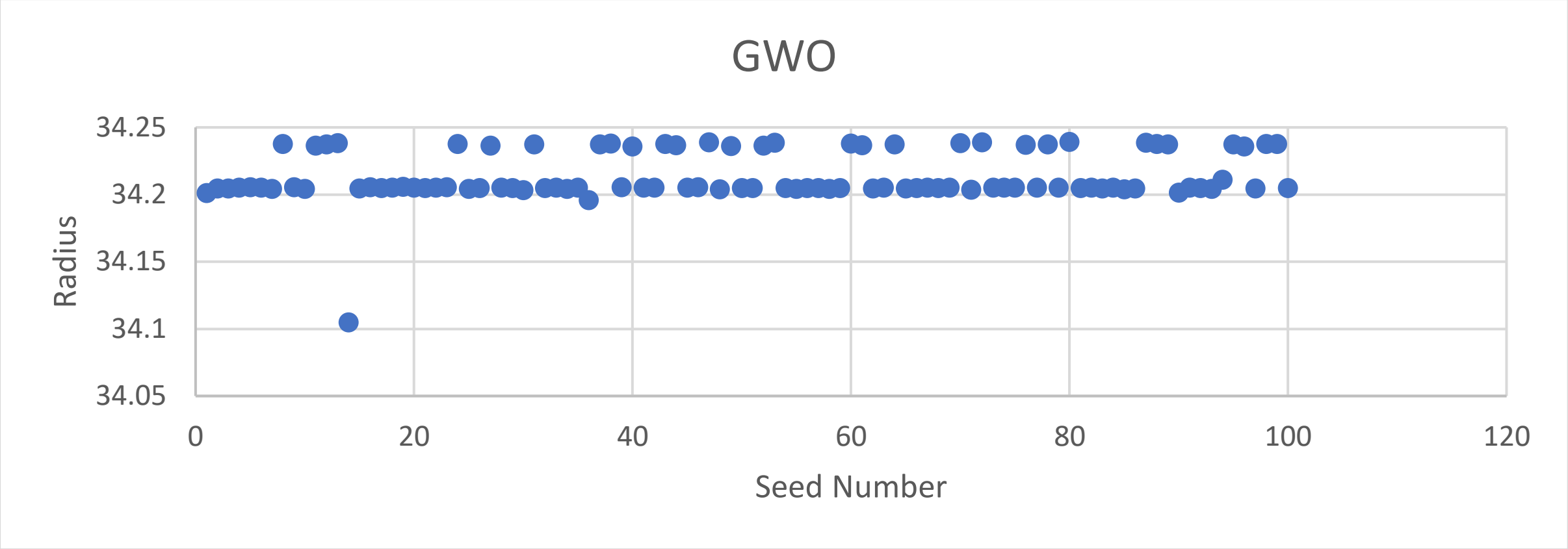

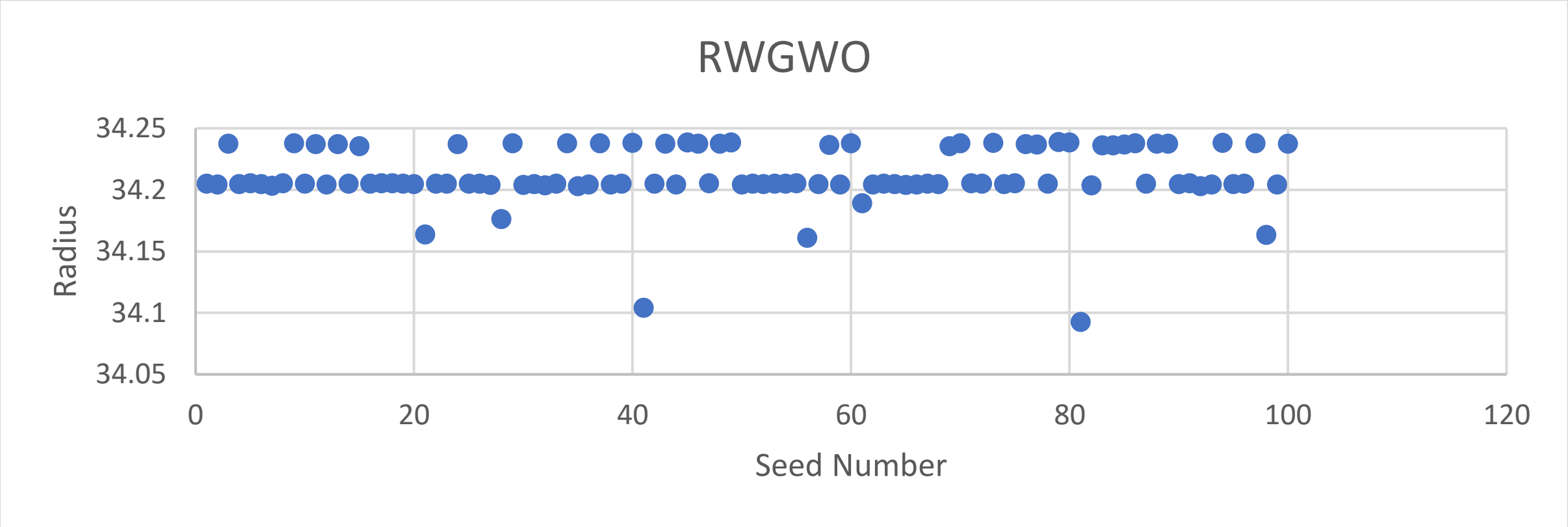

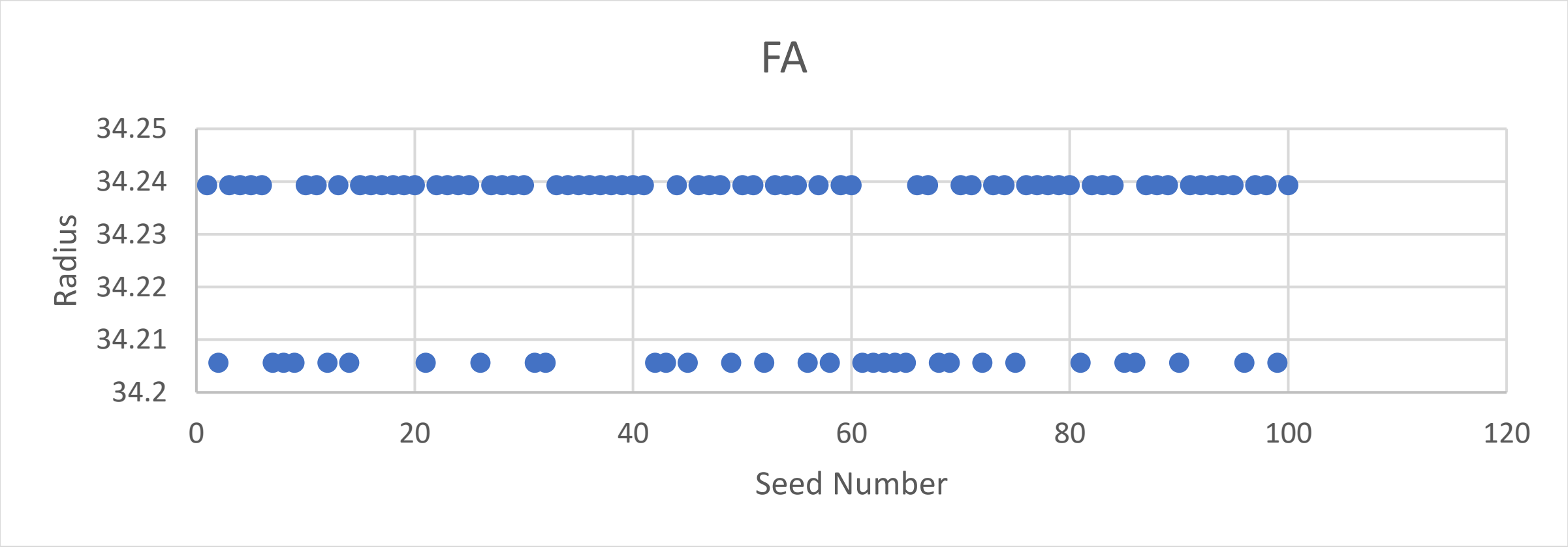

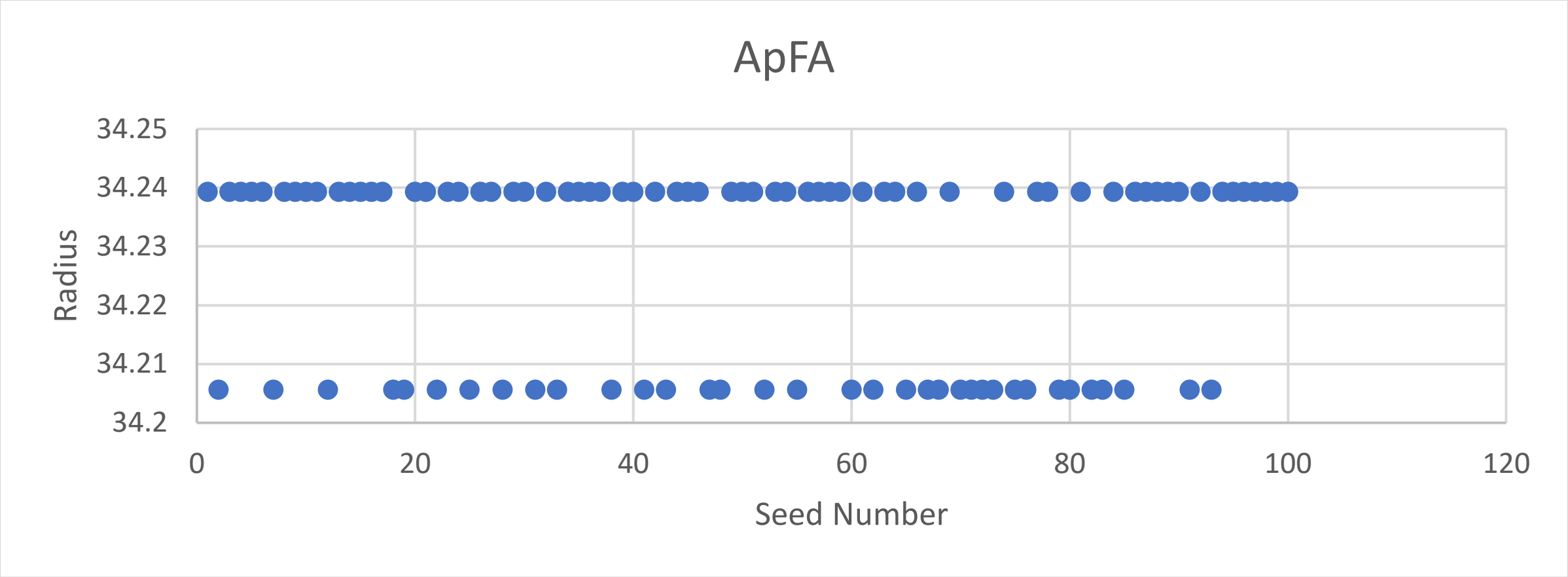

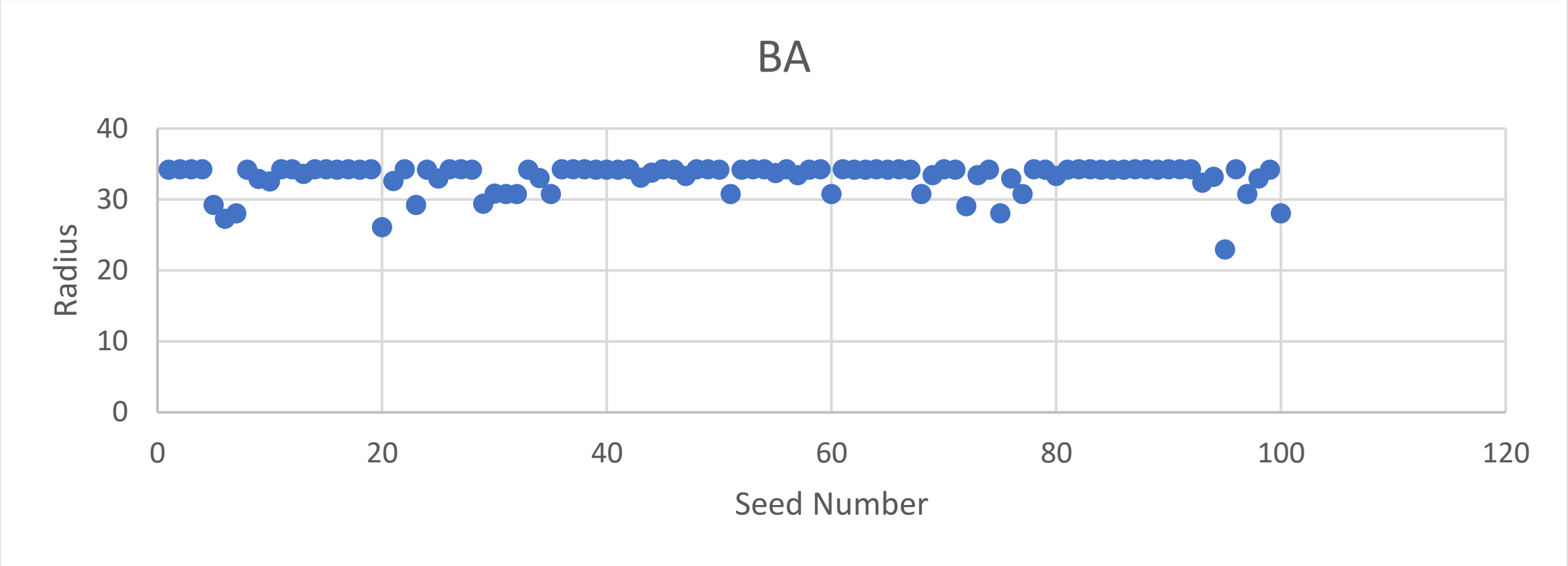

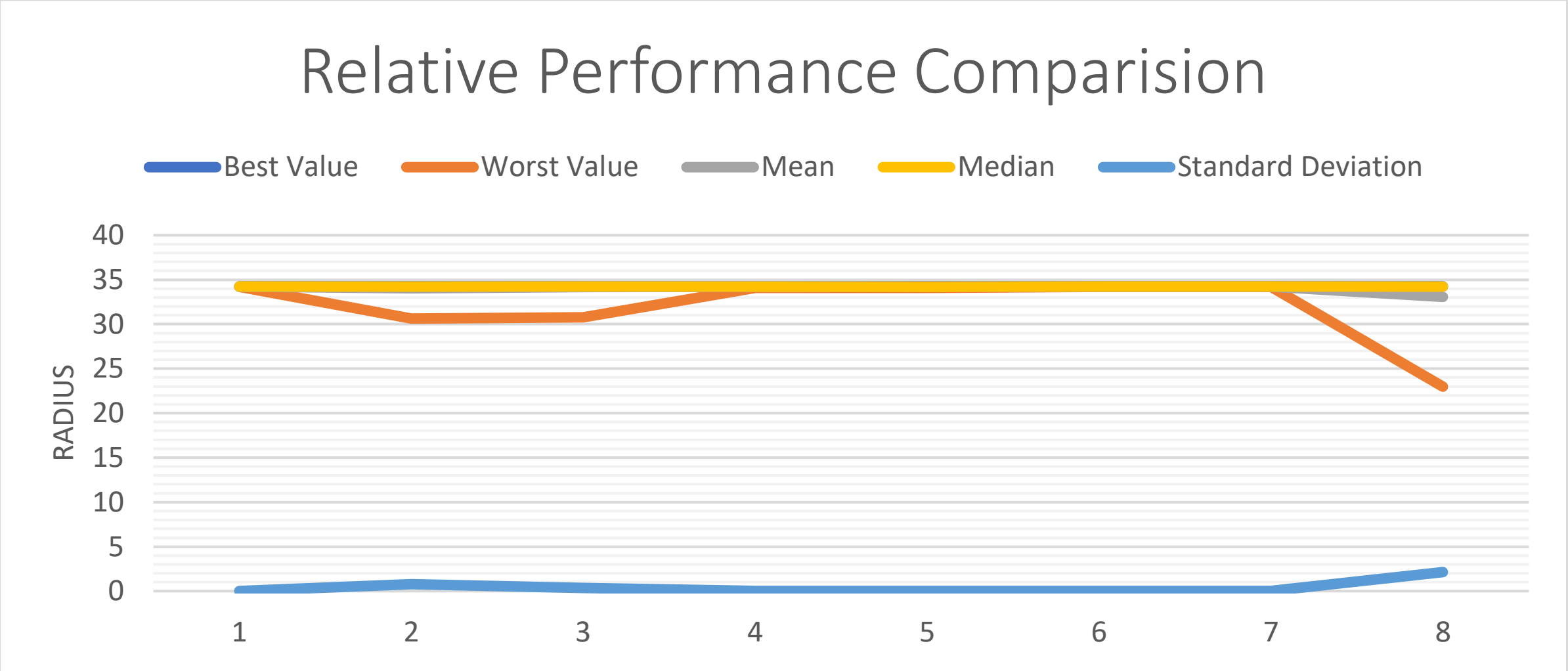

PSO, PSOd, CPSO, FA, ApFA, and BA have achieved the best value for the radius. BA achieves the worst value. In terms of average and standard deviation, there is a very minimal difference between ApFA, FA, and PSO, but PSO has been the best-fit algorithm for solving circle packing problems with these parameters. An interesting result to note here is that along with PSO, FA and ApFA median value of CPSO is also equal to the best value showing these algorithms' stability. Lastly, here also with the most considerable value for standard deviation BAT algorithm has depicted quiet instability.

Table13: Max Iterations = 100 & Number of Particles = 100

| | PSO | PSOd | CPSO | GWO | RWGWO | FA | ApFA | BA |
|---|---|---|---|---|---|---|---|---|
| Max Iterations | 100 | 100 | 100 | 100 | 100 | 100 | 100 | 100 |
| Number of Particles | 100 | 100 | 100 | 100 | 100 | 100 | 100 | 100 |
| Best Value | **34.2393** | **34.2393** | **34.2393** | 34.239 | 34.2382 | 24.3964 | **34.2393** | **34.2393** |
| Worst Value | 34.2056 | 33.7033 | 34.2056 | 33.945 | 34.0017 | **17.7777** | 34.2056 | 28.0431 |
| Mean | 34.2278 | 34.2163 | 34.2245 | 34.20204 | 34.2066 | 20.57748 | **34.22885** | 33.83874 |
| Median | **34.2393** | 34.2376 | **34.2393** | 34.20305 | 34.2032 | 20.3871 | **34.2393** | 34.2056 |
| Standard Deviation | 0.01596 | 0.05738 | 0.01681 | 0.043934 | 0.028737 | 1.246118 | **0.015665** | 1.159995 |

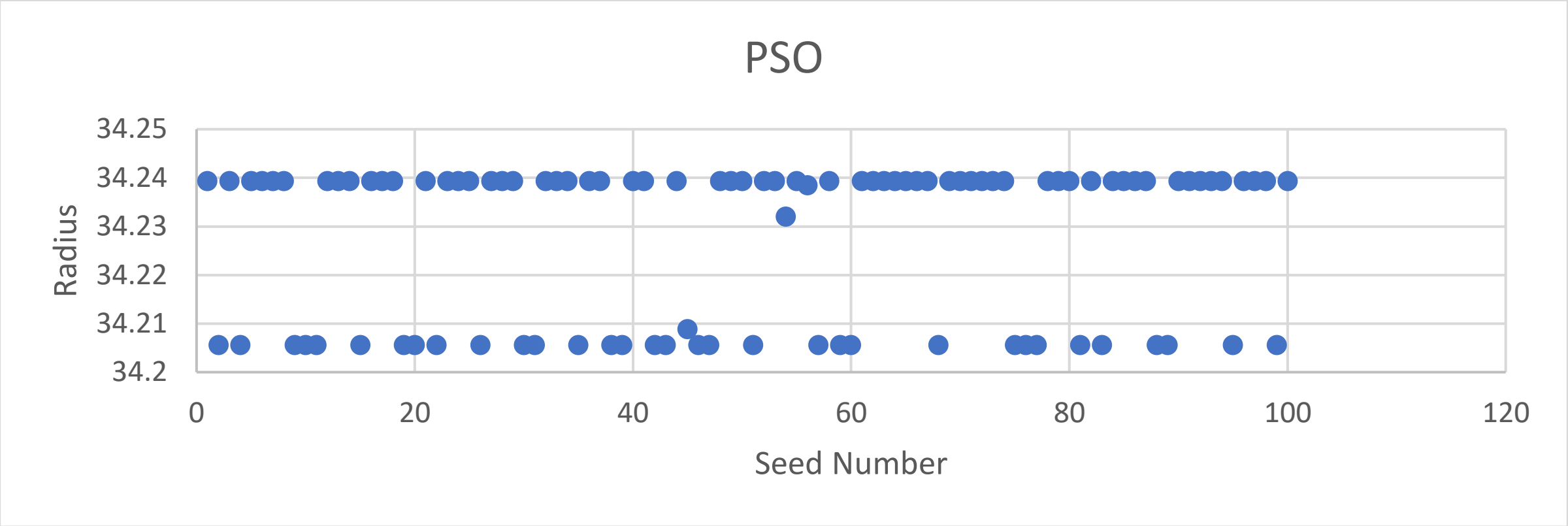

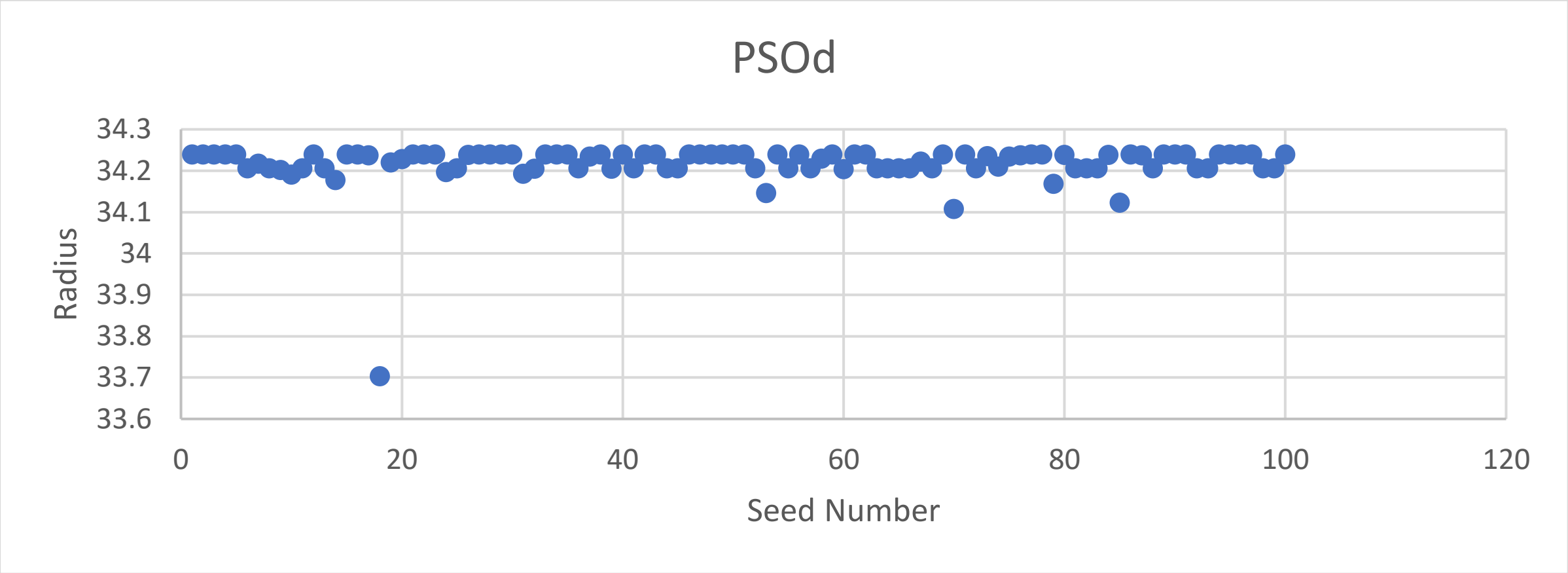

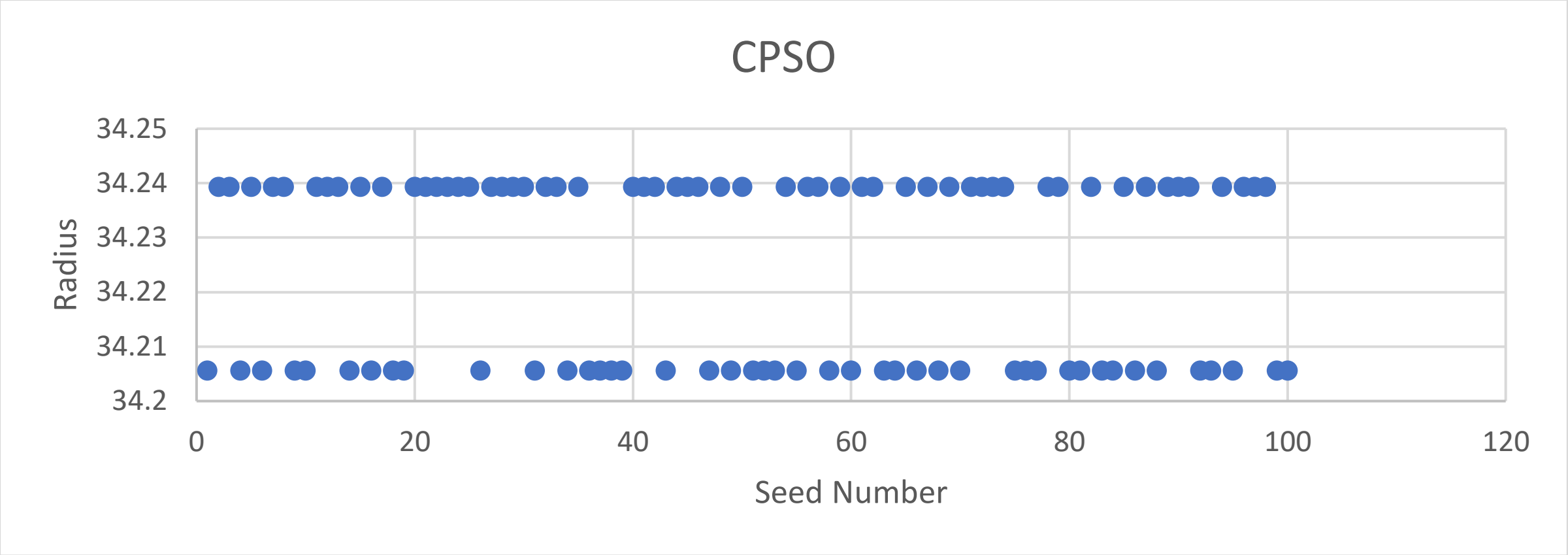
CPSO
34.25
34.24
34.23
34.22
34.21
34.2
Radius
0
20
40
60
80
100
120
Seed Number

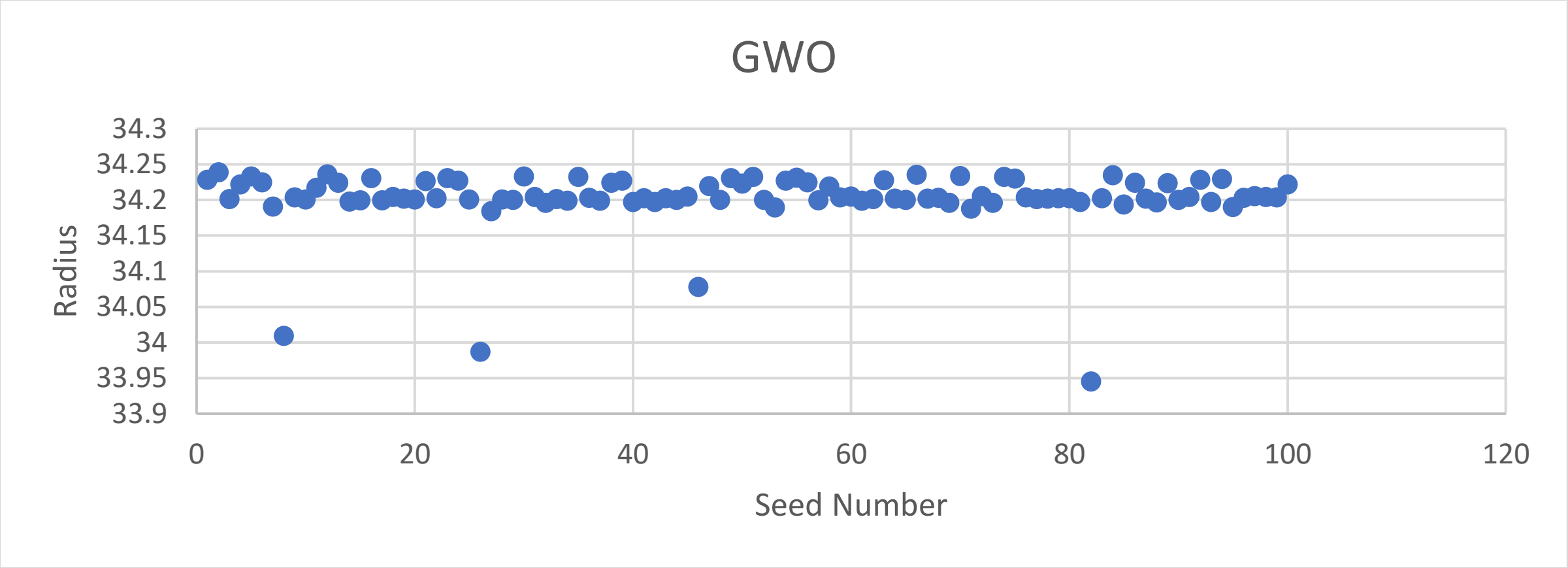
GWO
34.3
34.25
34.2
34.15
34.1
34.05
34
33.95
33.9
Radius
0
20
40
60
80
100
120
Seed Number

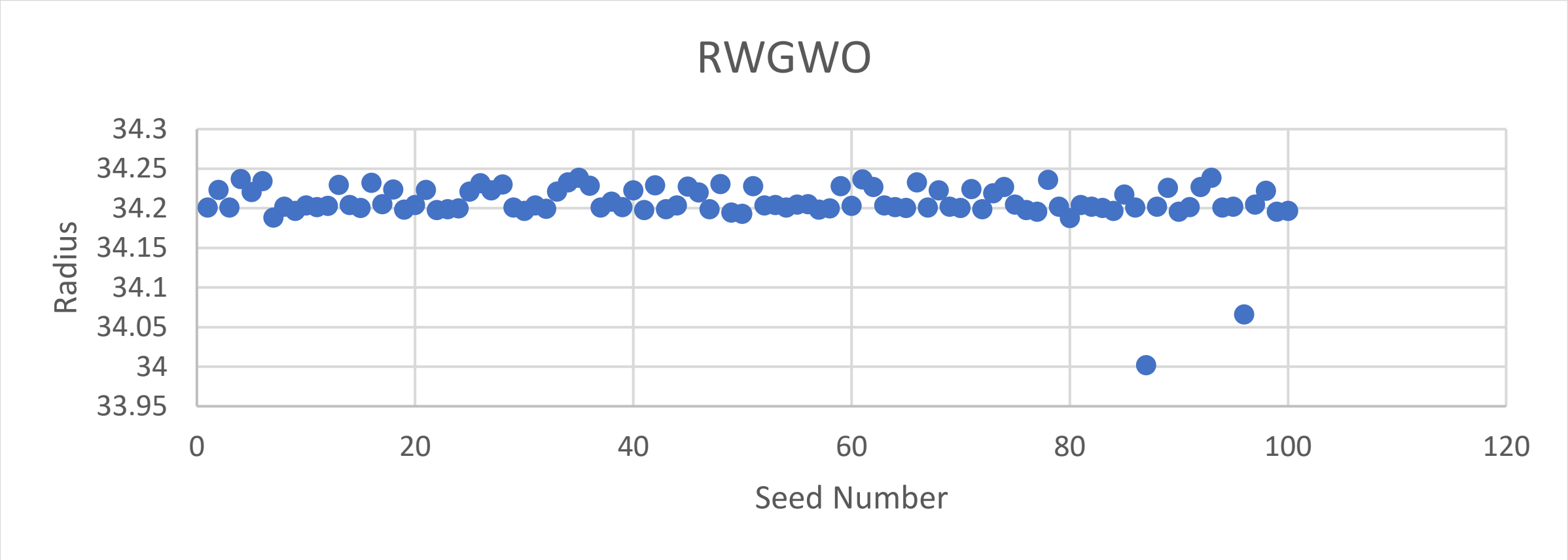
RWGWO
34.3
34.25
34.2
34.15
34.1
34.05
34
33.95
Radius
0
20
40
60
80
100
120
Seed Number

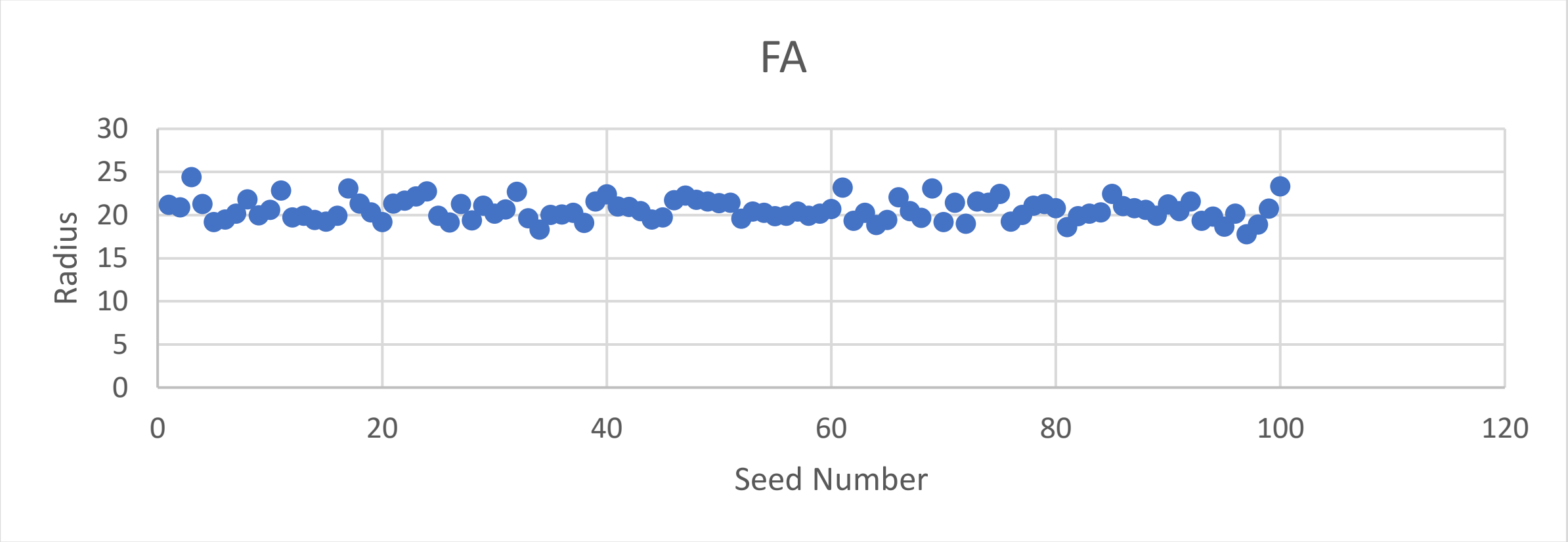
FA
30
25
20
15
10
5
0
Radius
0
20
40
60
80
100
120
Seed Number

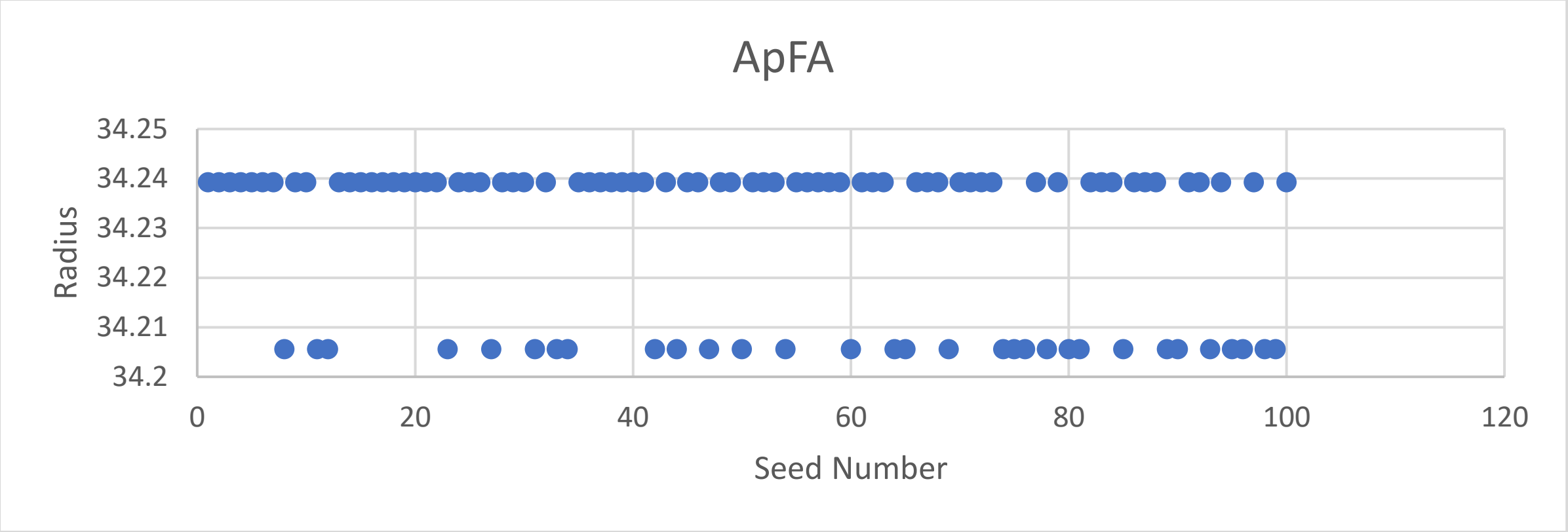

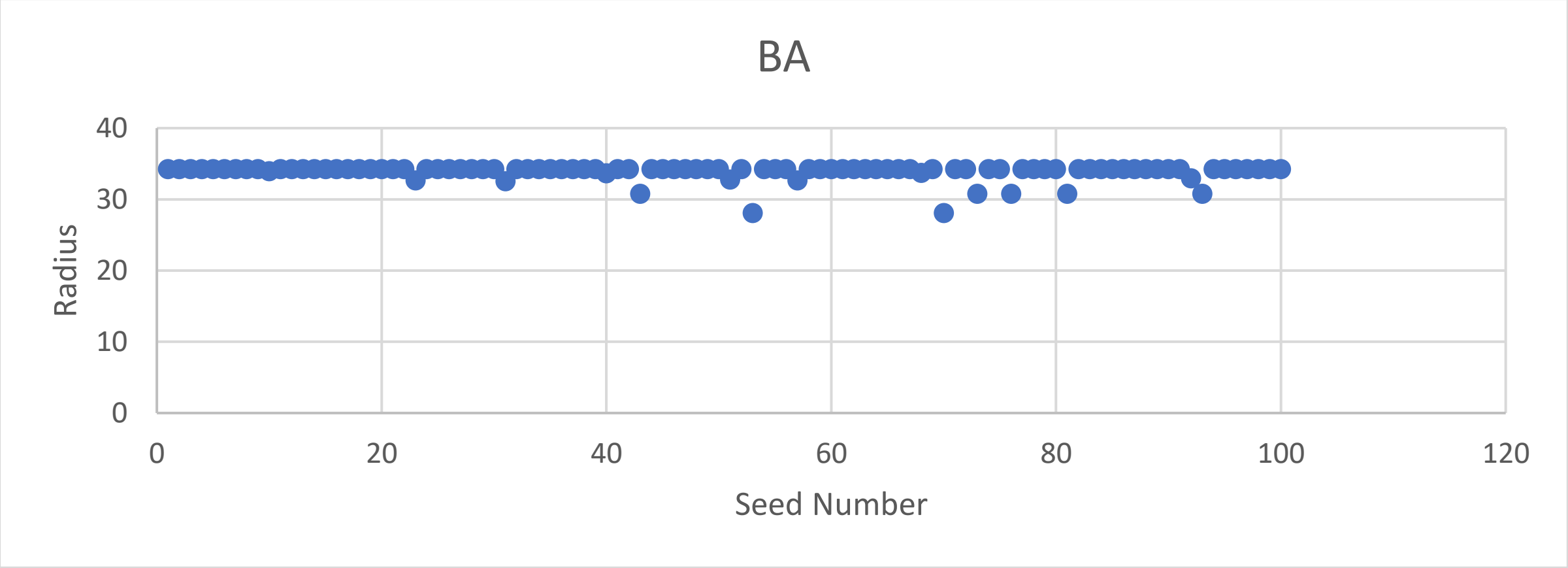

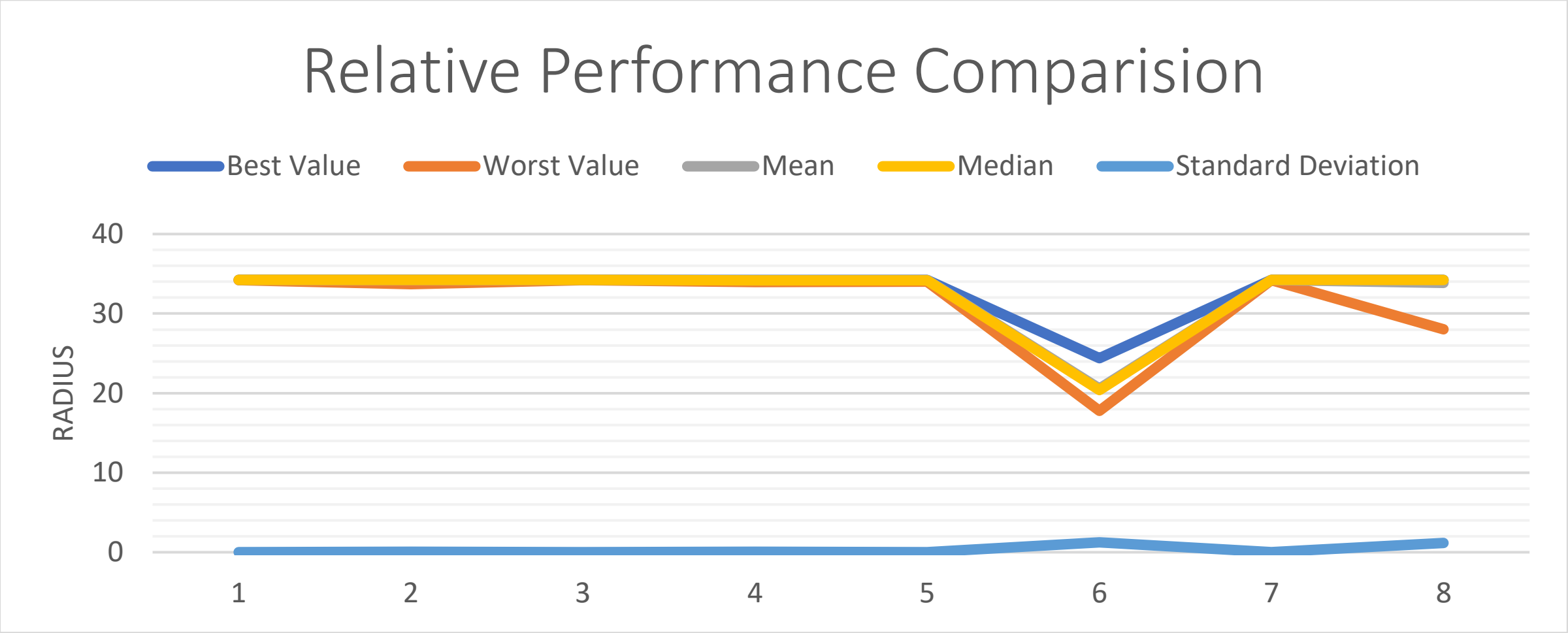

PSO, PSOd, CPSO, ApFA, and BA have achieved the best value for the radius. FA reaches the worst value. In terms of average and standard deviation, there is a very minimal difference between ApFA, CPSO, and PSO. Still, ApFA has been the best-fit algorithm for solving the circle packing problem with these parameters. Notably, with the increase in the number of particles, the performance of BA has been improved to a quiet extent and has achieved some stability. Due to slower convergence, FA has shown poor results.

Table14: Max Iterations = 500 & Number of Particles = 100

| | PSO | PSOd | CPSO | GWO | RWGWO | FA | ApFA | BA |
|---|---|---|---|---|---|---|---|---|
| Max Iterations | 500 | 500 | 500 | 500 | 500 | 500 | 500 | 500 |

| Number of Particles | 100 | 100 | 100 | 100 | 100 | 100 | 100 | 100 |
|---|---|---|---|---|---|---|---|---|
| Best Value | **34.2393** | **34.2393** | **34.2393** | 34.239 | 34.2389 | 34.2392 | **34.2393** | **34.2393** |
| Worst Value | 34.2056 | 30.7809 | 34.2056 | 34.0491 | 34.0937 | 34.2055 | 34.2056 | **28.0432** |
| Mean | **34.237** | 34.186 | 34.2258 | 34.20955 | 34.21439 | 34.23179 | 34.2329 | 33.83235 |
| Median | **34.2393** | **34.2285** | **34.2393** | **34.205** | **34.20515** | **34.2392** | **34.2393** | **34.2056** |
| Standard Deviation | **0.00834** | 0.34491 | 0.01659 | 0.028235 | 0.02403 | 0.014027 | 0.013287 | 1.15153 |

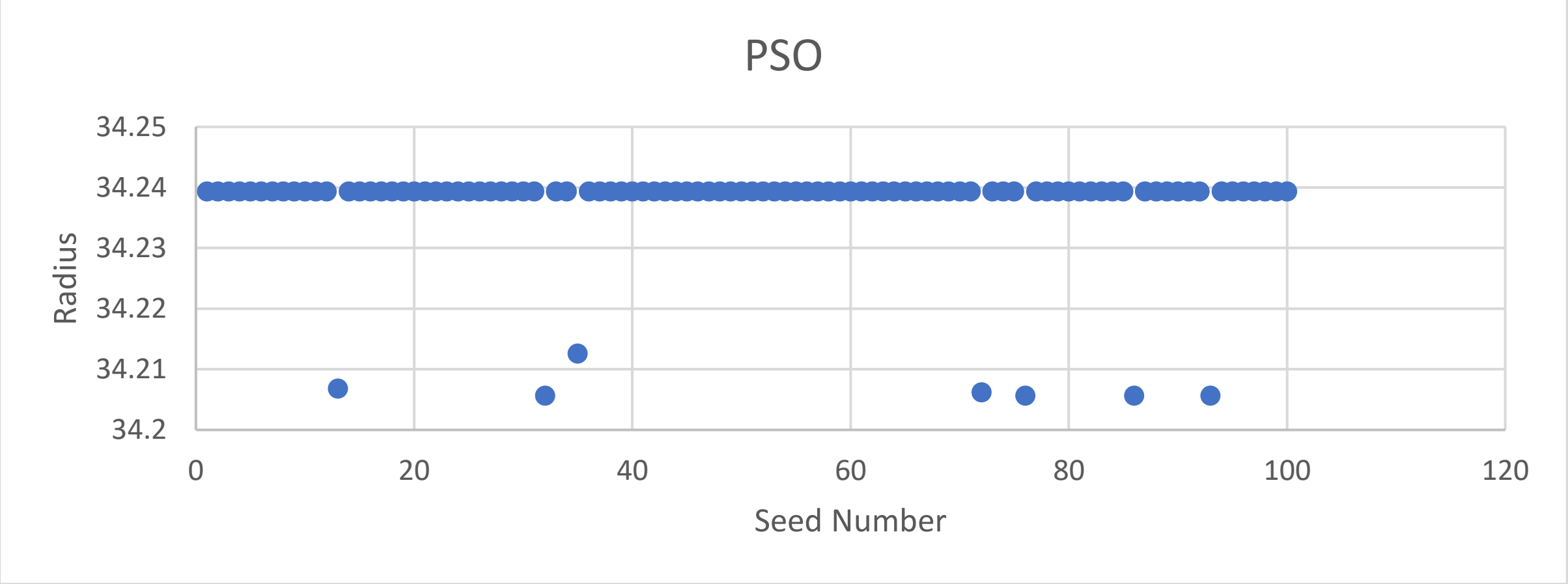

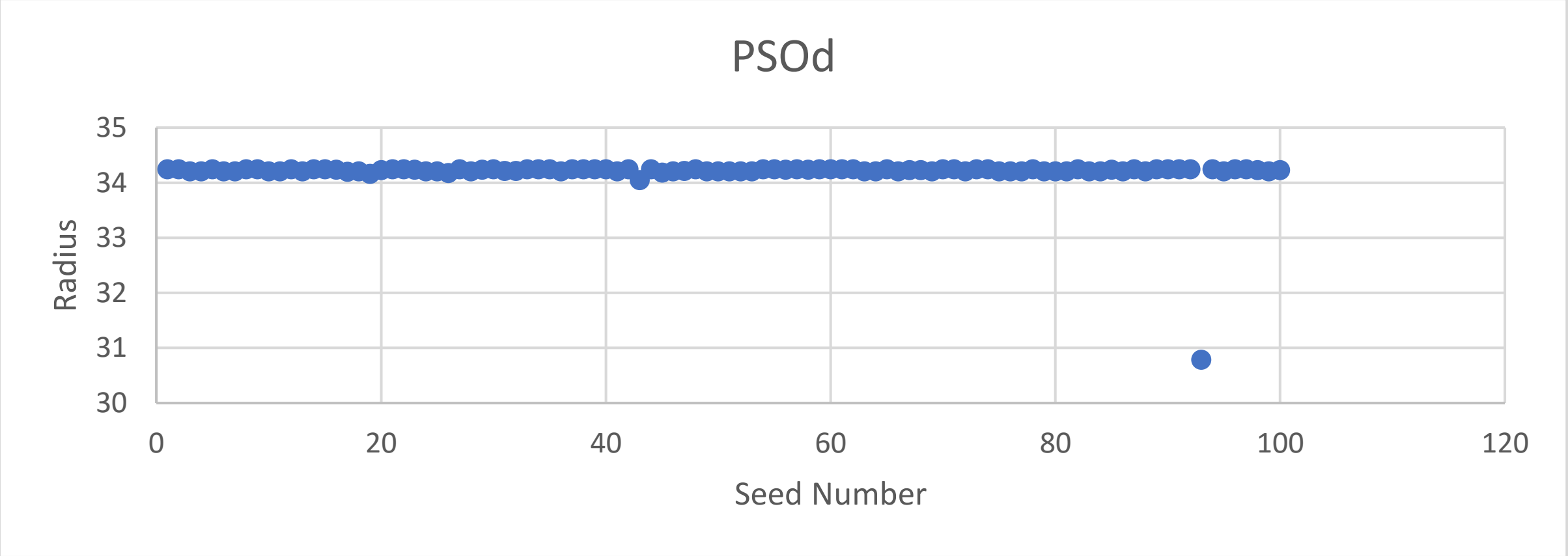

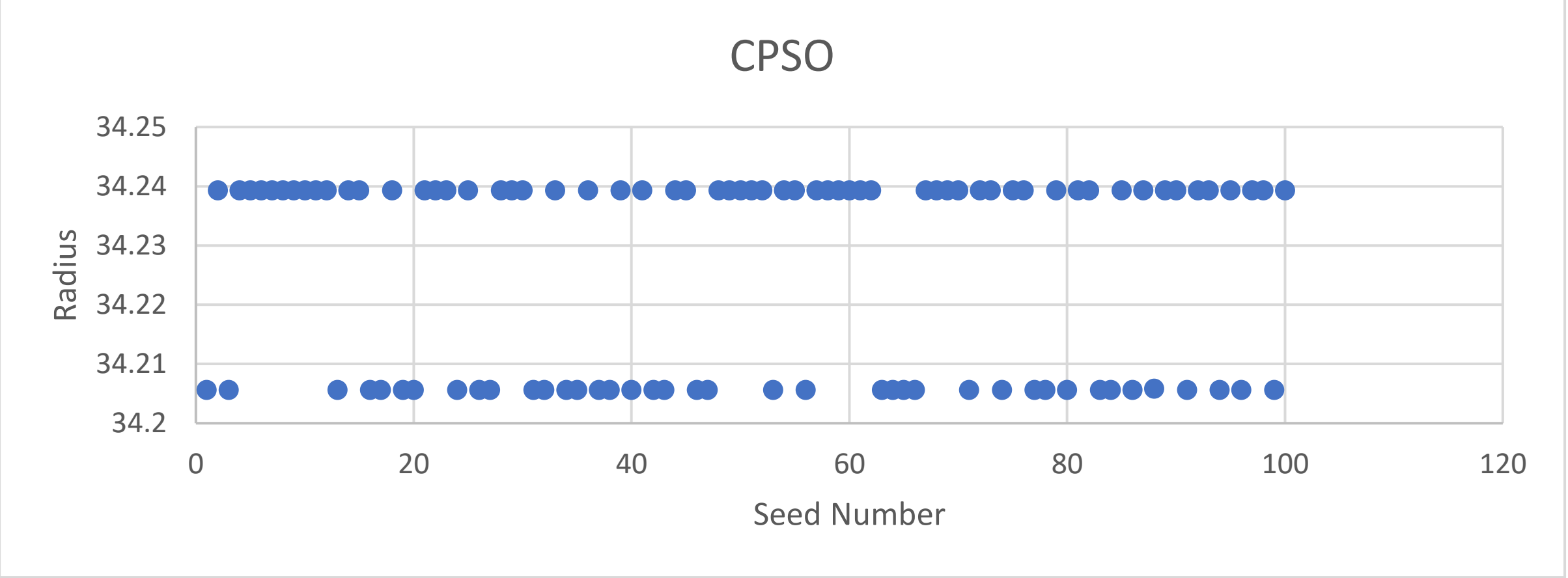

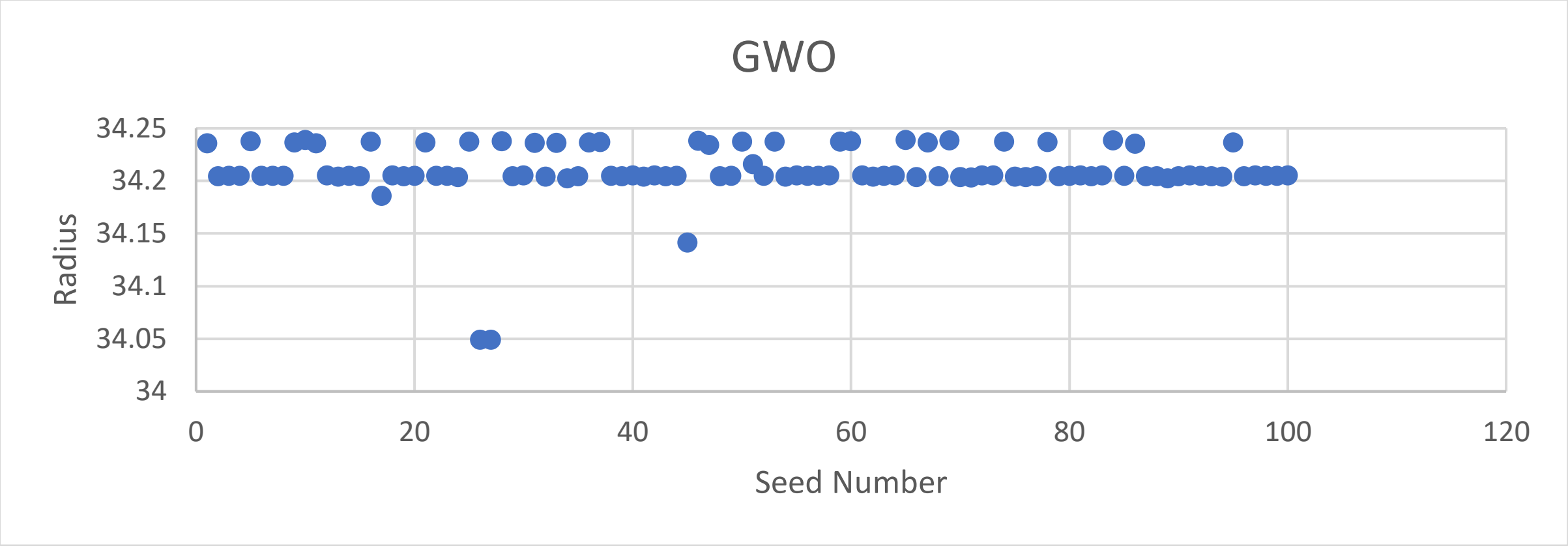
GWO
34.25
34.2
34.15
34.1
34.05
34
Radius
0
20
40
60
80
100
120
Seed Number

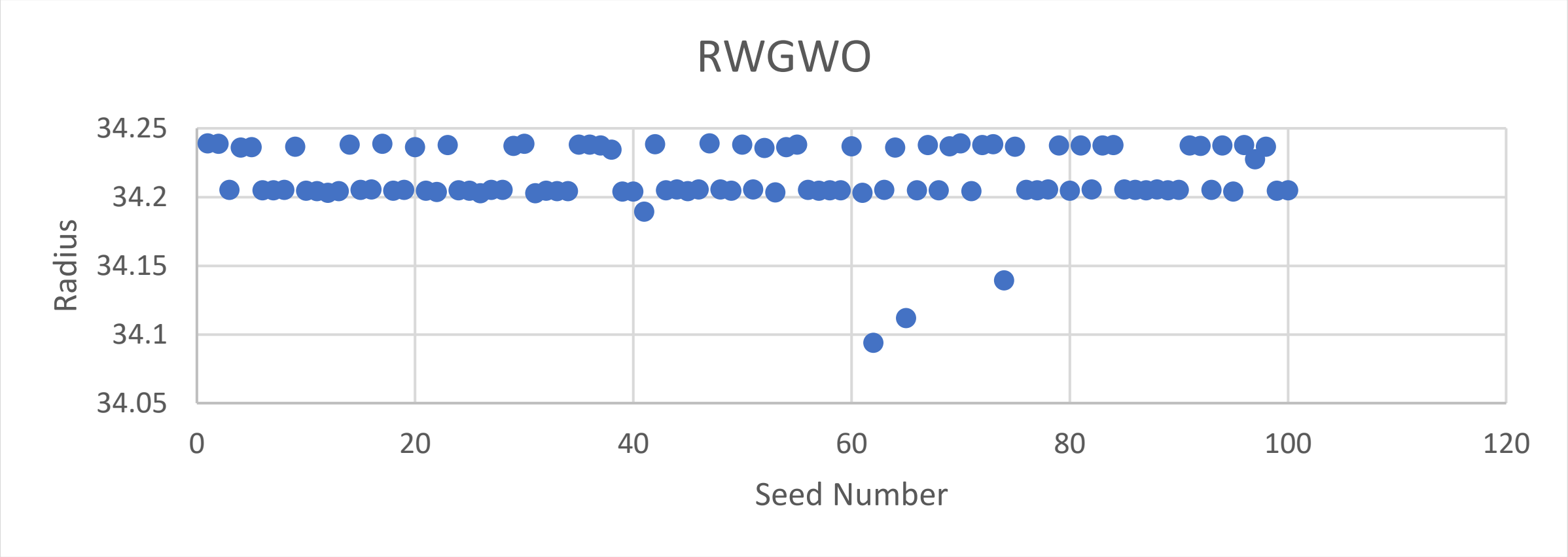
RWGWO
34.25
34.2
34.15
34.1
34.05
Radius
0
20
40
60
80
100
120
Seed Number

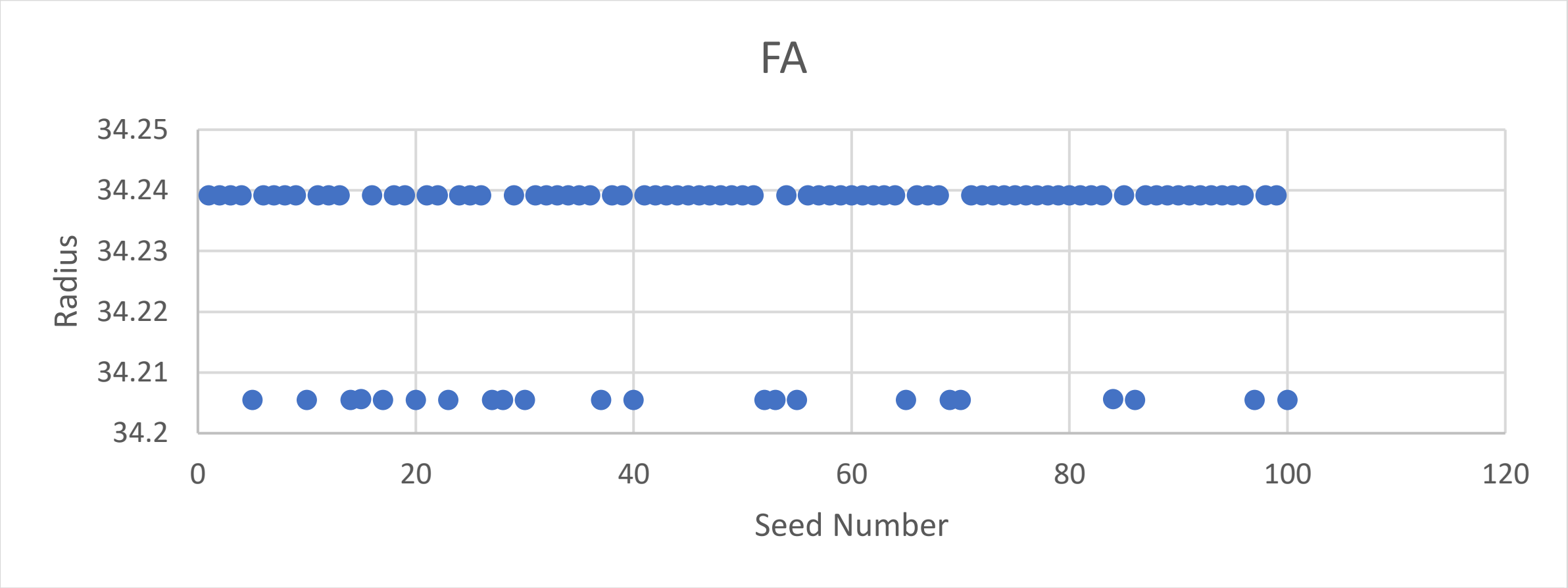
FA
34.25
34.24
34.23
34.22
34.21
34.2
Radius
0
20
40
60
80
100
120
Seed Number

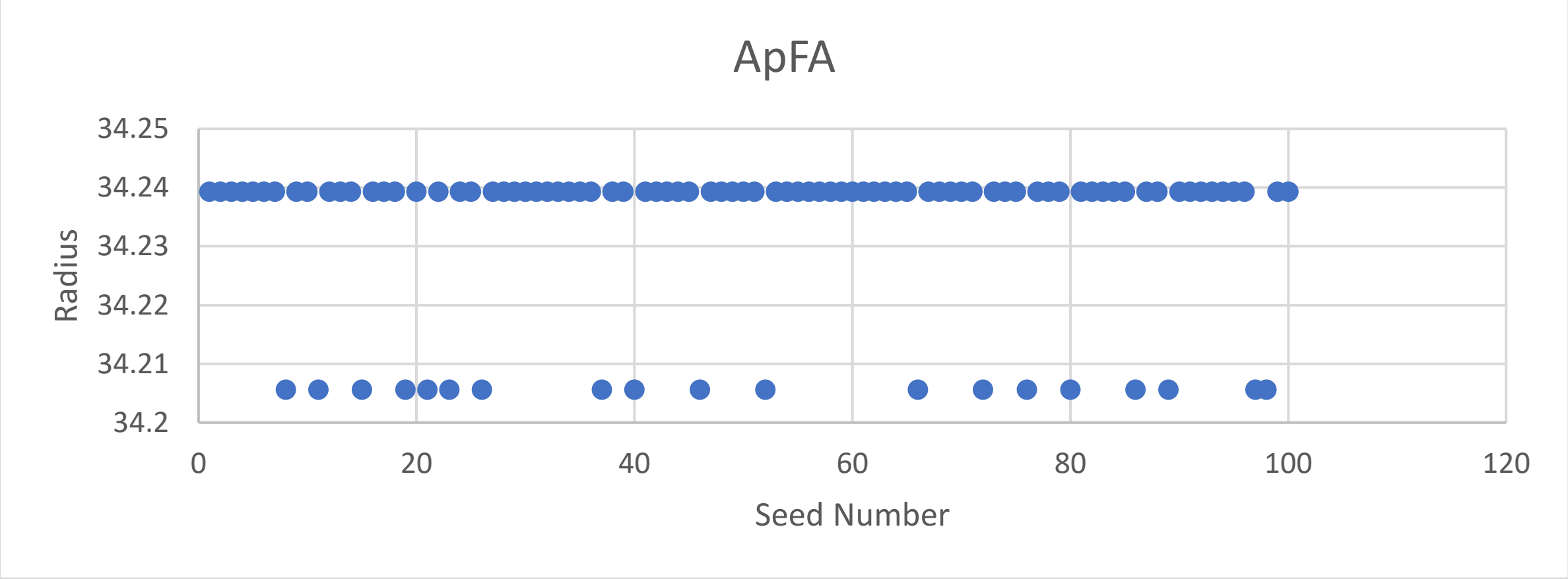

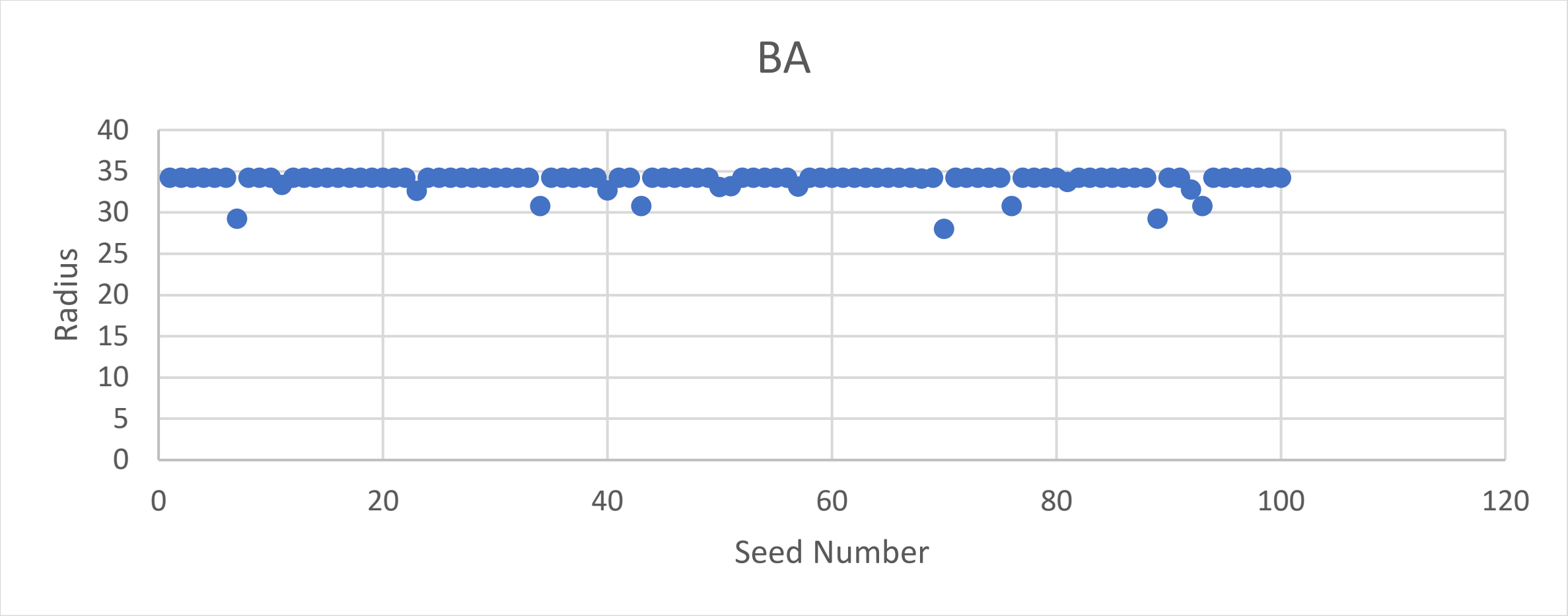

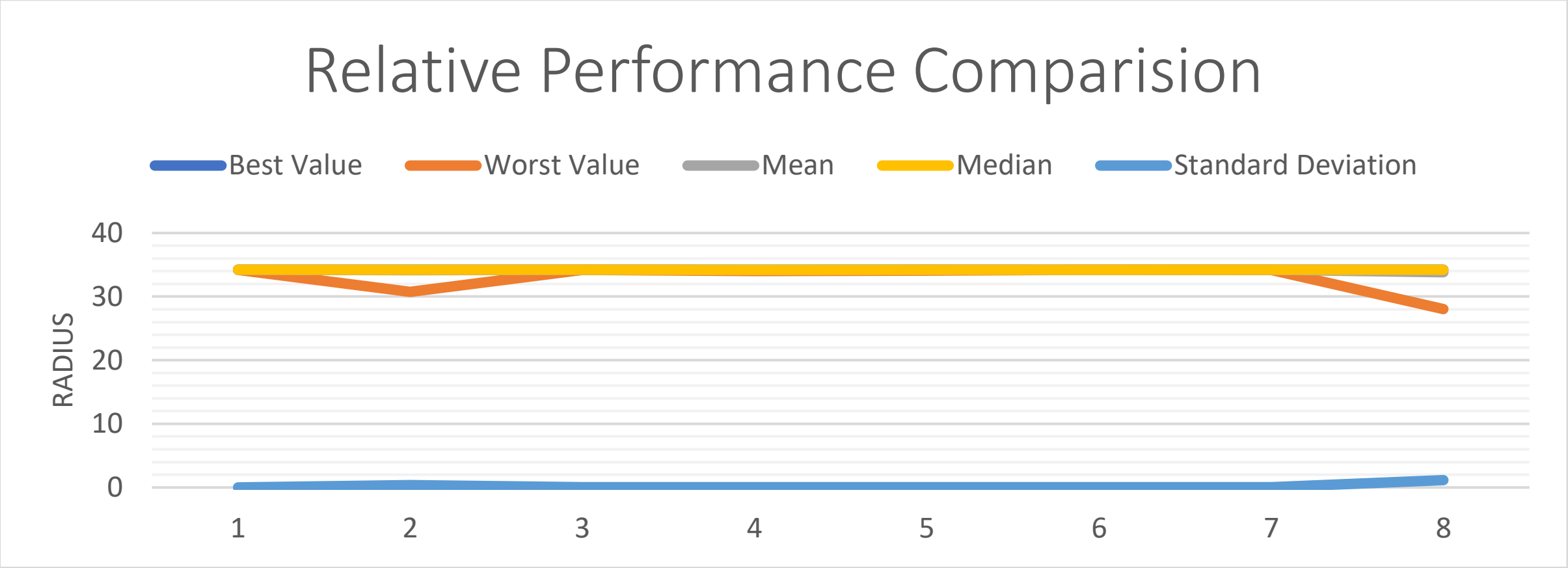

PSO, PSOd, CPSO, ApFA, and BA have achieved the best value for the radius. GWO, RWGWO, and FA have also reached the best value with an almost negligible difference. BA achieves the worst value. In terms of average and standard deviation, PSO has been the best-fit algorithm for solving circle packing problems with these parameters.

**Table15** Max Iterations = 1000 & Number of Particles = 100

| | PSO | PSOd | CPSO | GWO | RWGWO | FA | ApFA | BA |
|---|---|---|---|---|---|---|---|---|
| Max Iterations | 1000 | 1000 | 1000 | 1000 | 1000 | 1000 | 1000 | 1000 |
| Number of Particles | 100 | 100 | 100 | 100 | 100 | 100 | 100 | 100 |
| Best Value | **34.2393** | **34.2393** | **34.2393** | 34.2392 | 34.2392 | **34.2393** | **34.2393** | **34.2393** |
| Worst Value | 34.2056 | 30.7809 | 34.2056 | 34.1317 | 34.1347 | 34.2056 | 34.2056 | **28.0432** |
| Mean | **34.2383** | 34.1878 | 34.2299 | 34.21623 | 34.21992 | 34.23391 | 34.23391 | 33.97153 |
| Median | **34.2393** | **34.2352** | **34.2393** | **34.2054** | **34.2055** | **34.2393** | **34.2393** | **34.2393** |
| Standard Deviation | **0.00578** | 0.34516 | 0.0152 | 0.018731 | 0.018594 | 0.012417 | 0.012417 | 0.894556 |

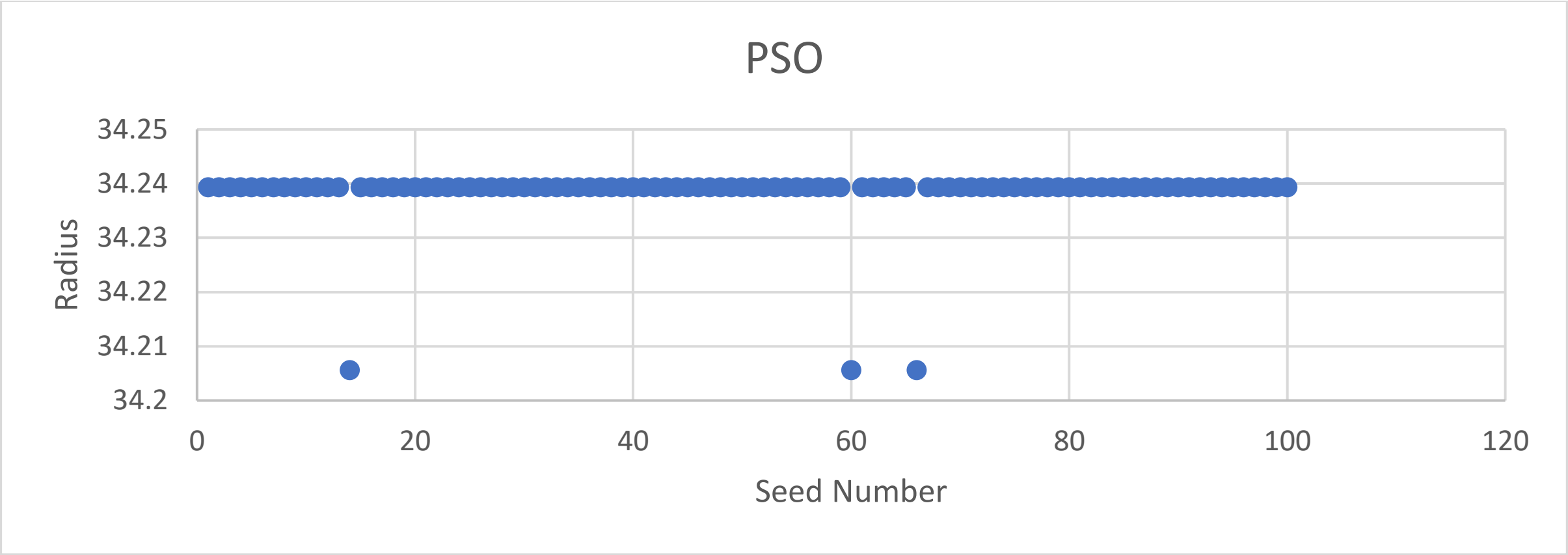
PSO
34.25
34.24
34.23
34.22
34.21
34.2
Radius
0
20
40
60
80
100
120
Seed Number

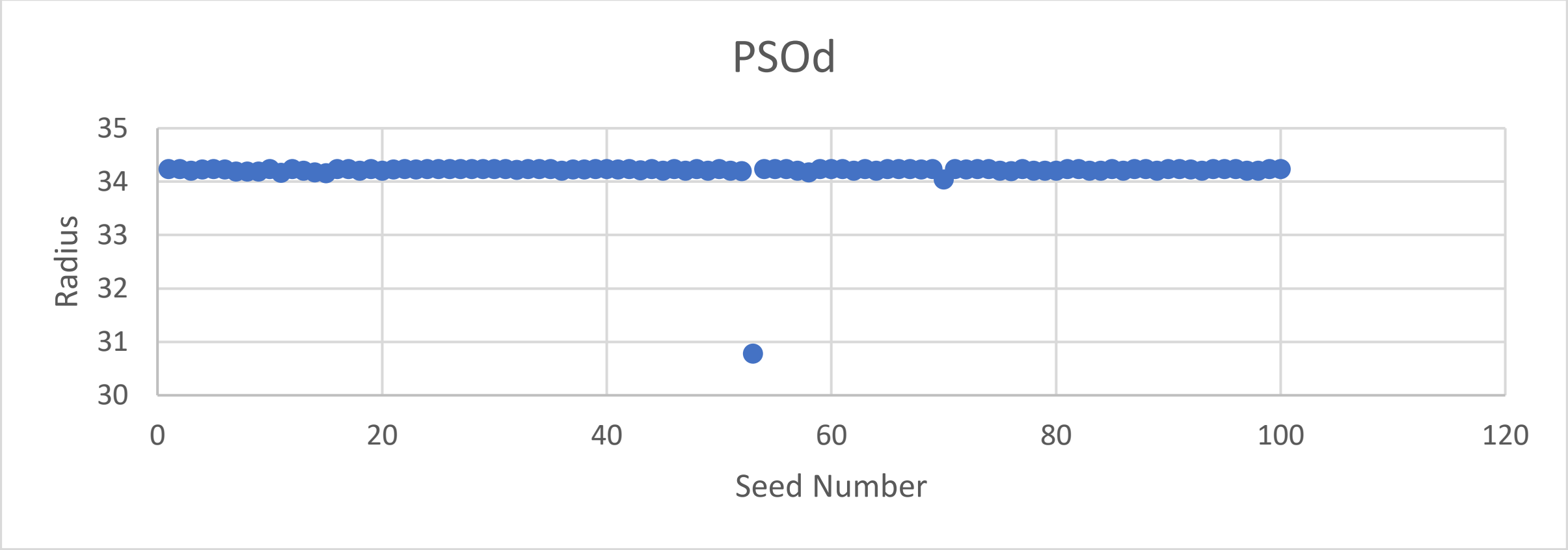
PSOd
35
34
33
32
31
30
Radius
0
20
40
60
80
100
120
Seed Number

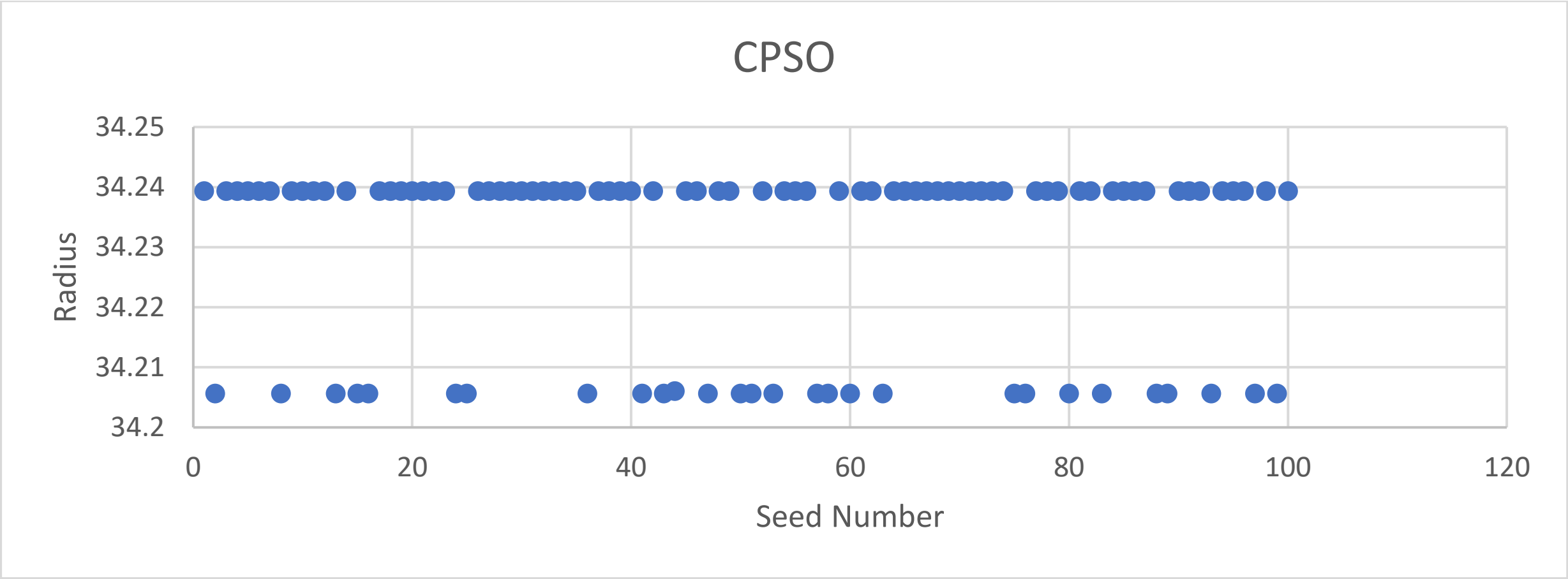
CPSO
34.25
34.24
34.23
34.22
34.21
34.2
Radius
0
20
40
60
80
100
120
Seed Number

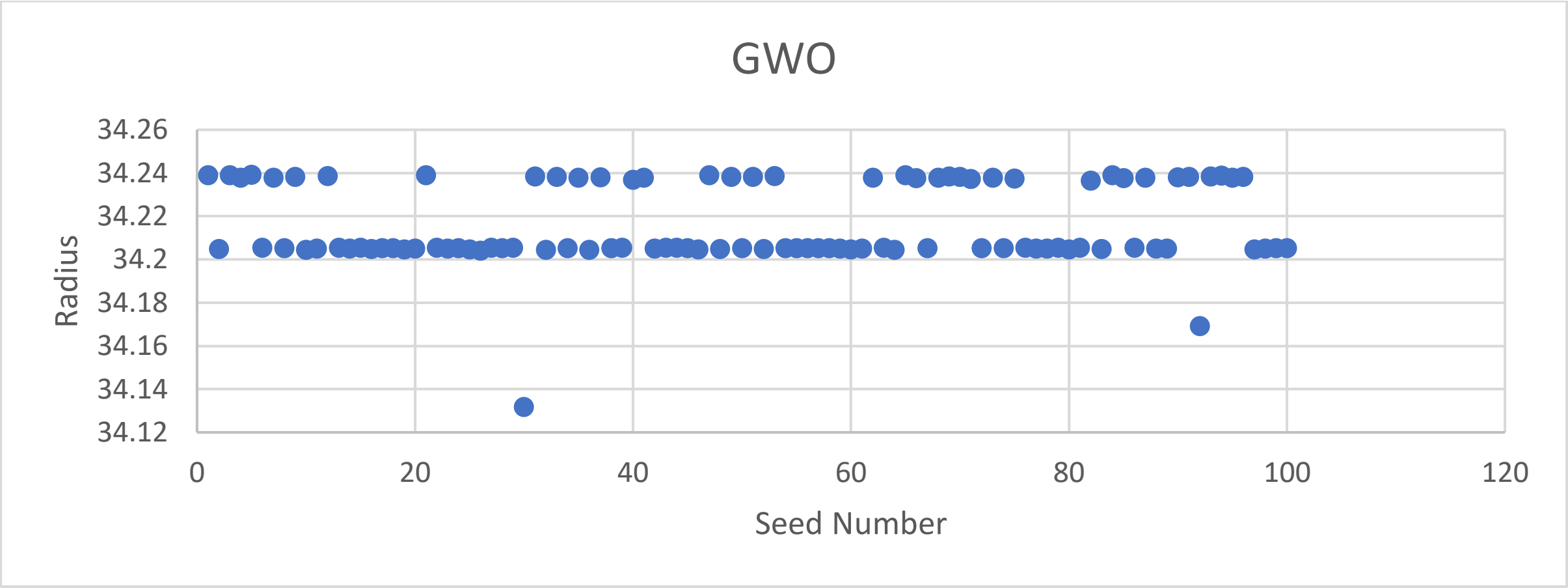
GWO
34.26
34.24
34.22
34.2
34.18
34.16
34.14
34.12
Radius
0
20
40
60
80
100
120
Seed Number

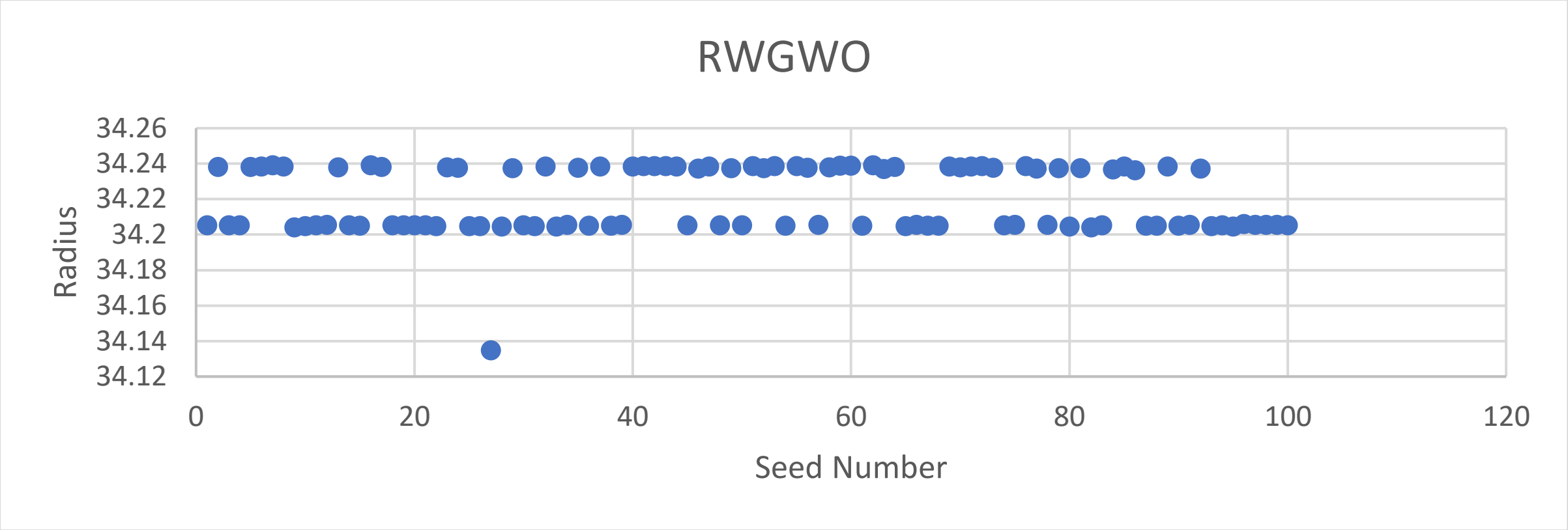
RWGWO
34.26
34.24
34.22
34.2
34.18
34.16
34.14
34.12
Radius
0
20
40
60
80
100
120
Seed Number

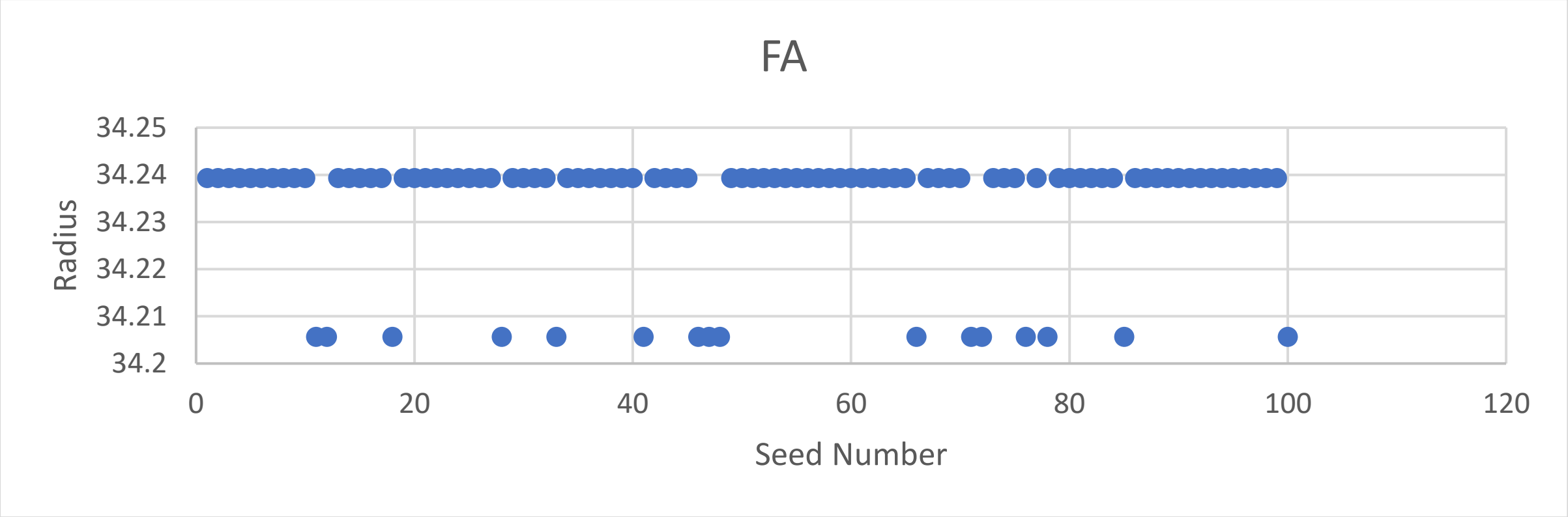
FA
34.25
34.24
34.23
34.22
34.21
34.2
Radius
0
20
40
60
80
100
120
Seed Number

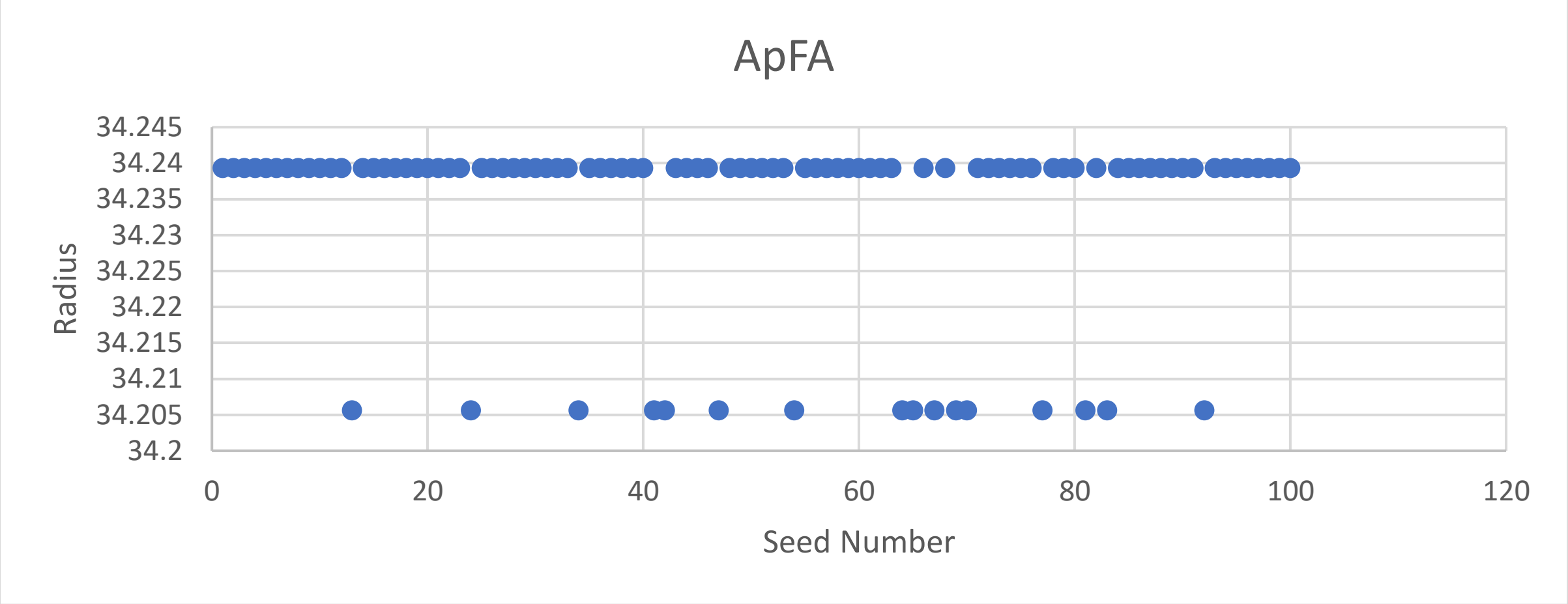
ApFA
34.245
34.24
34.235
34.23
34.225
34.22
34.215
34.21
34.205
34.2
Radius
0
20
40
60
80
100
120
Seed Number

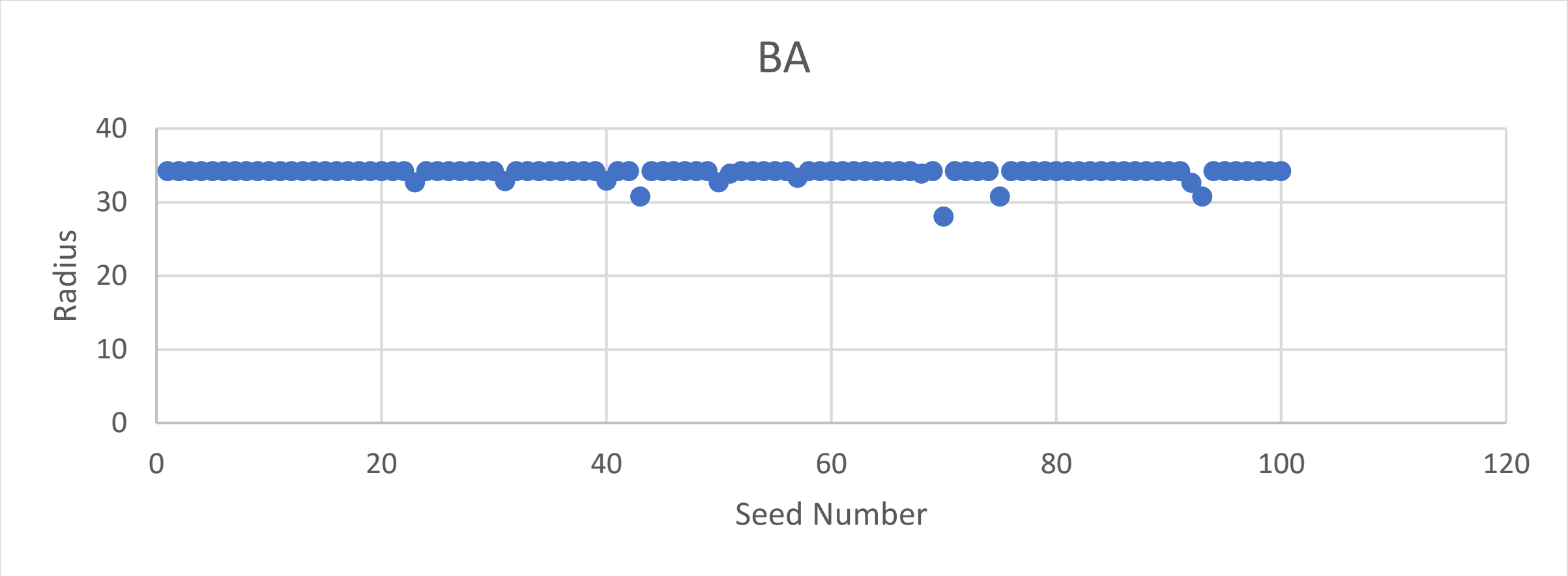

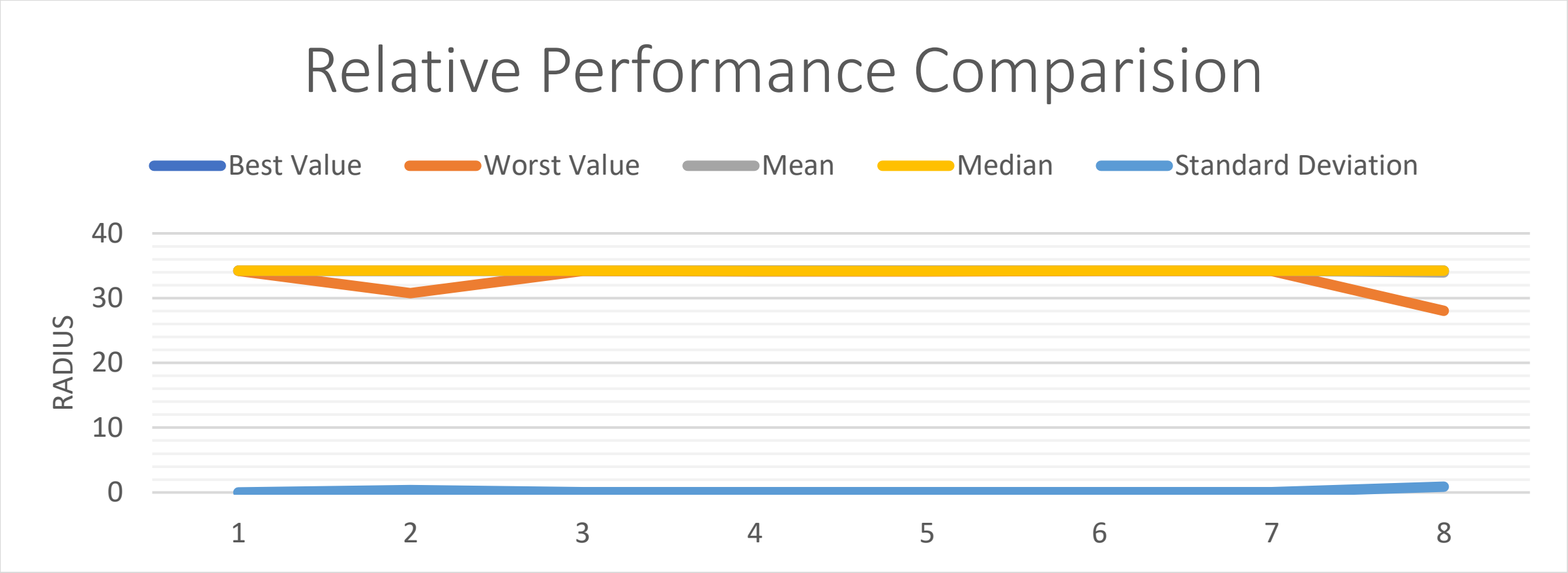

PSO, PSOd, CPSO, ApFA, and BA have achieved the best value for the radius. GWO, RWGWO, and FA have also achieved the best value with an almost negligible difference. BA reaches the worst value. Among all the six combinations, PSO has achieved the highest average and lowest deviation value with these parameters making it the best fit for solving the circle packing problem. In addition to median values of PSO, CPSO, ApFA, and FA, interestingly, the median value of BA also equals maximum value, which means with an increase in the number of particles and maximum iterations, the BA algorithm becomes more stable reliable.

**6. Conclusion and Future Works**

The best fit and worst fit algorithms for each parameter set are summarized in the table below. Although PSO has shown relatively good results overall, we have seen that no one such algorithm can be regarded as best fit for all parameter sets, showcasing the problem's highly non-linear and sensitive nature. PSO can be considered the most promising algorithm for solving this problem, but it depends on the data set and parameters given the susceptible nature of the problem.

**Table16** Best and Worst algorithms

| S. No. | No. of particles | Max. no of iterations | Best algorithm | Second best algorithm | Worst algorithm |
|---|---|---|---|---|---|
| 1 | 50 | 100 | PSO | ApFA | FA |
| 2 | 50 | 500 | FA | ApFA | BA |
| 3 | 50 | 1000 | PSO | FA | BA |
| 4 | 100 | 100 | ApFA | PSO | FA |
| 5 | 100 | 500 | PSO | ApFA | BA |
| 6 | 100 | 1000 | PSO | FA, ApFA | BA |

Since each algorithm was implemented on 100 different seeds, efficacy is also calculated as the number of times the algorithm achieves the best value (i.e., the maximum radius- 34.239) under the given set of Maximum Iteration and Number of Particles. Efficacy for best and 2$^{nd}$ best algorithms is depicted below in graphs. It can be concluded that the efficacy in general increases with an increase in the number of particles and number of iterations; however, the trade-off here is computational time and resources.

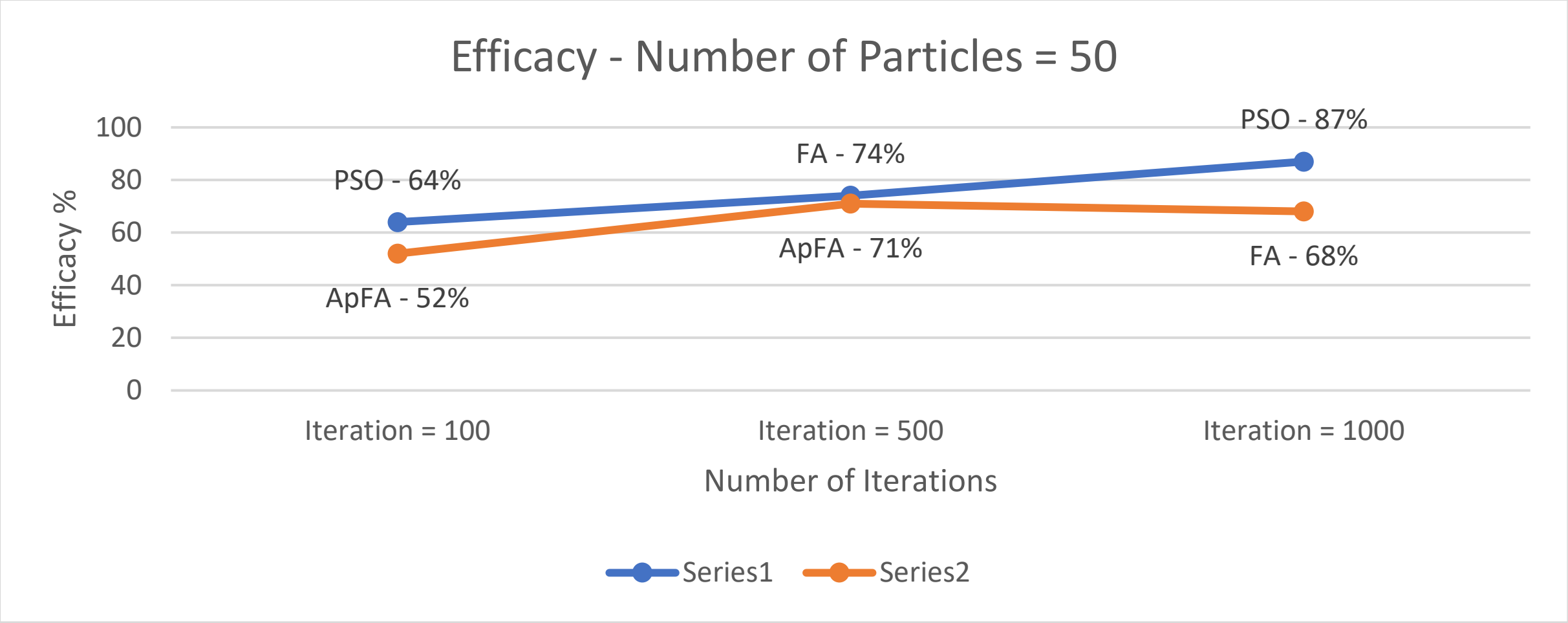

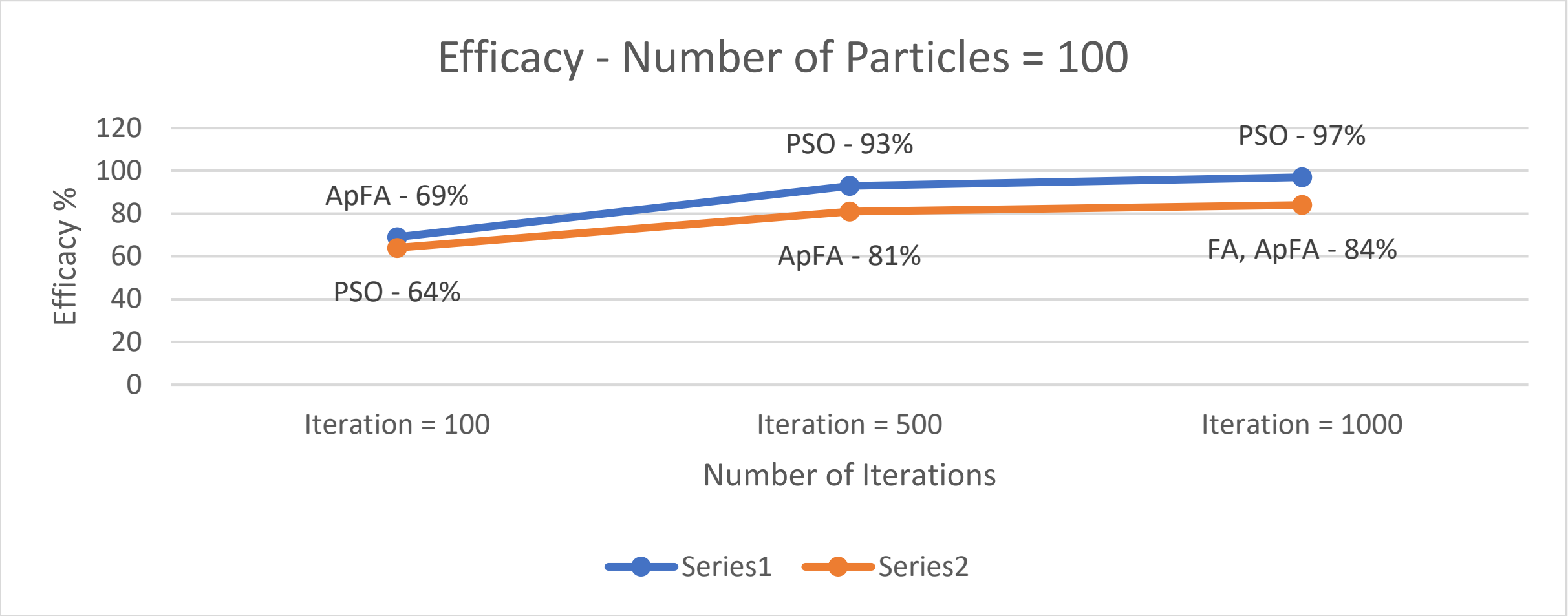

Circle packing problems have received a considerable amount of attention in the mathematics literature, but only modest attention in the operational research literature. This is surprising, considering the numerous areas of existing and potential applications: Radiation Treatment Planning, Container Loading, Cylinder Packing, Social Distancing Problems, Tree Plantation etc. We have tried to implement Meta-Heuristic algorithm along with some recent advancement on a test data set and the future works include implementation of these approaches on an industrial (realistic) circle packing problem.